\documentclass{article}%
\usepackage{amsfonts}
\usepackage{amsmath}
\usepackage{amssymb}
\usepackage{graphicx}%
\providecommand{\U}[1]{\protect\rule{.1in}{.1in}}

\begin{document}

\title{ Bernstein Diffusions for a Class of Linear Parabolic Partial Differential Equations}
\author{Pierre-A. Vuillermot$^{\ast}$ and Jean-C. Zambrini$^{\ast\ast}$\\UMR-CNRS 7502, Institut \'{E}lie Cartan, Nancy, France$^{\ast}$\\Grupo de F\'{\i}sica Matem\'{a}tica da Universidade de Lisboa, Portugal$^{\ast
\ast}$}
\date{}
\maketitle

\begin{abstract}
In this article we prove the existence of Bernstein processes which we
associate in a natural way with a class of non-autonomous linear parabolic
initial- and final-boundary value problems defined in bounded convex subsets
of Euclidean space of arbitrary dimension. Under certain conditions regarding
their joint endpoint distributions, we also prove that such processes become
reversible Markov diffusions. Furthermore we show that those diffusions
satisfy two It\^{o} equations for some suitably constructed Wiener processes,
and from that analysis derive Feynman-Kac representations for the solutions to
the given equations. We then illustrate some of our results by considering the
heat equation with Neumann boundary conditions both in a one-dimensional
bounded interval and in a two-dimensional disk.

\end{abstract}

\section{Introduction and Outline}

It is well known that It\^{o}'s theory of stochastic differential equations
makes it possible to associate a Markov process to essentially any
second-order elliptic differential operator, and that the fine properties of
each one of these objects allow one to get precise information about the
other. In particular, the knowledge of the behavior of such a diffusion
typically leads to the discovery of new phenomena regarding the solutions to a
host of elliptic and parabolic partial differential equations, ranging from
Dirichlet and Neumann initial-boundary value problems to equations describing
wave front propagation in periodic and random media. This is testified, for
instance, by the many results and references in \cite{dynkin} and
\cite{freidlin}.

By the same token it is also possible to associate Markov processes with the
Schr\"{o}dinger equation of quantum physics in a variety of ways, as in
\cite{carlen} where the author's considerations rest on the principles of
stochastic mechanics set forth in \cite{nelson1}, and in \cite{zambrini} whose
constructions are related to the properties of the stochastic variational
calculus introduced in \cite{yasue}.

The diffusion processes constructed in \cite{carlen} and \cite{zambrini} are,
however, quite different from the more traditional dissipative diffusions in
that they encode all the conservative properties of quantum mechanics. In
particular, they exhibit the invariance under time reversal inherent in the
Schr\"{o}dinger equation and thereby maintain a perfect symmetry between past
and future. More to the point, the Markov processes of \cite{zambrini}
actually emerge as particular cases of reversible diffusions that belong to
the larger class of the so-called reciprocal or Bernstein processes, whose
theory was launched many years ago in \cite{bernstein} following
Schr\"{o}dinger's seminal contribution in \cite{schroedinger}. The theory of
Bernstein processes was subsequently further developed and systematically
investigated in \cite{jamison}, and since then has played an important
r\^{o}le in relating various fields such as the Malliavin calculus and
Euclidean quantum mechanics, or Markov bridges with jumps and L\'{e}vy
processes, to name only a few (see for instance \cite{cruzwuzambrin},
\cite{cruzeirozambrini}, \cite{gulicasteren}, \cite{privaultzambrin},
\cite{vancasteren} and the references therein for a more complete account).

Given these facts there remains the interesting and intriguing question
whether it is possible to associate reversible diffusion processes to
parabolic equations of general form whose solutions typically display
irreversible behavior, and to get new and nontrivial information out of this
association. In his attempt to understand certain analogies between the
properties of Brownian motion and quantum mechanics in the last section of
\cite{schroedinger}, Schr\"{o}dinger answers the question positively by
analyzing a simple case through statistical arguments.

It is the purpose of this article to show that this can also be achieved by
purely analytical means and indeed with a considerable degree of generality.
Let $D\subset\mathbb{R}^{d}$ be a bounded open convex subset whose boundary is
denoted by $\partial D$. We consider parabolic \textit{initial-boundary value
problems }of the form%
\begin{align}
\partial_{t}u(x,t)  &  =\frac{1}{2}\operatorname{div}_{x}\left(
k(x,t)\nabla_{x}u(x,t)\right)  -\left(  l(x,t),\nabla_{x}u(x,t)\right)
_{\mathbb{R}^{d}}-V(x,t)u(x,t),\nonumber\\
\text{\ }(x,t)  &  \in D\times\left(  0,T\right]  ,\nonumber\\
u(x,0)  &  =\varphi(x),\text{ \ \ }x\in D,\nonumber\\
\frac{\partial u(x,t)}{\partial n_{k}(x,t)}  &  =0,\text{ \ \ \ }%
(x,t)\in\partial D\times\left(  0,T\right]  \label{parabolicproblem}%
\end{align}
with $T\in\left(  0,+\infty\right)  $ arbitrary, $(.,.)_{\mathbb{R}^{d}}$ the
Euclidean inner product in $\mathbb{R}^{d}$, and where the last relation in
(\ref{parabolicproblem}) stands for the conormal derivative of $u$ relative to
the matrix-valued function $k$. Furthermore $l$, $V$ and $\varphi$ are an
$\mathbb{R}^{d}$-valued vector-field and real-valued functions, respectively.

Let us now associate with (\ref{parabolicproblem}) its \textit{adjoint
final-boundary value problem, }namely,%
\begin{align}
-\partial_{t}v(x,t)  &  =\frac{1}{2}\operatorname{div}_{x}\left(
k(x,t)\nabla_{x}v(x,t)\right)  +\operatorname{div}_{x}\left(
v(x,t)l(x,t)\right)  -V(x,t)v(x,t),\nonumber\\
(x,t)  &  \in D\times\left[  0,T\right)  ,\nonumber\\
v(x,T)  &  =\psi(x),\text{ \ \ }x\in D,\nonumber\\
\frac{\partial v(x,t)}{\partial n_{k}(x,t)}  &  =0,\text{ \ \ \ }%
(x,t)\in\partial D\times\left[  0,T\ \right)  , \label{adjointproblem}%
\end{align}
in order to build reversibility and eventually the Markov property into the
theory we develop below.

We then organize the remaining part of this article in the following way: in
Section 1 we prove the existence of Bernstein processes $Z_{\tau\in\left[
0,T\right]  }$ wandering in $\overline{D}:=D\cup\partial D$, which we
associate with (\ref{parabolicproblem}) and (\ref{adjointproblem}) in a very
natural way provided the coefficients therein be sufficiently smooth. We also
show there that under certain conditions which pertain to their joint initial
and final distributions, the processes in question become reversible Markov
diffusions. Moreover, we can express their probability density as the product
of the solutions to (\ref{parabolicproblem}) and (\ref{adjointproblem}), and
prove that they also satisfy two It\^{o} stochastic differential equations for
some suitably constructed Wiener processes. This last fact is not merely
anecdotic, as it allows us eventually to get Feynman-Kac representations of
those solutions in terms of $Z_{\tau\in\left[  0,T\right]  }$ in a perfectly
symmetric manner, by invoking It\^{o}'s \textit{backward} stochastic calculus
in the case of (\ref{parabolicproblem}), and quite independently the
usual\textit{ forward} stochastic calculus in the case of
(\ref{adjointproblem}). In Section 1 we also explain why the notion of
reversible Markov diffusion which emerges from our considerations corresponds
to a generalization of the classical notion of reversibility put forward in
\cite{kolmogorov} and inspired by the last section of \cite{schroedinger},
which was many years later reformulated in \cite{dobsukfritz} and
\cite{ikedawatanabe}. In Section 2 we illustrate some of our results by means
of two examples involving the simple heat equation in a one-dimensional
bounded interval and in a two-dimensional disk, and in one of them the process
we construct shares some of the features of a reflected diffusion. Finally,
for the sake of completeness and for the convenience of the reader, we devote
an appendix to reviewing the existence theory of weak solutions to
(\ref{parabolicproblem}) and (\ref{adjointproblem}) which provide the
ingredients that are essential to our construction when the coefficients
therein are indeed smooth enough.

Throughout this article we use the standard notations for all the functional
spaces we need without any further comments, referring the reader for instance
to \cite{adamsfournier}, \cite{eideliva}, \cite{eidelzhitara} and
\cite{tanabe}. We also use freely results on martingales and Wiener processes
from \cite{gihmanskohorod}.

\section{A Class of Bernstein Diffusions}

We are looking for processes $Z_{\tau\in\left[  0,T\right]  }$ whose natural
state space is the compact convex set $\overline{D}$ endowed with the Borel
$\sigma$-algebra $\mathcal{B(}\overline{D})$, and whose behavior corresponding
to time values belonging to any subinterval $\left(  s,t\right)
\subset\left[  0,T\right]  $ is conditioned by the knowledge of $Z_{s}$ and
$Z_{t}$ alone. This means that all past information gathered prior to time $s$
is irrelevant, as is all future information accumulated after time $t$. The
precise notion we need is the following (see \cite{jamison} and some of the
references therein for other equivalent formulations):

\bigskip

\textsc{Definition 1.} \textit{We say the }$\overline{D}$\textit{-valued
process }$Z_{\tau\in\left[  0,T\right]  }$ \textit{defined on the complete
probability space }$\left(  \Omega,\mathcal{F},\mathbb{P}\right)  $\textit{ is
a Bernstein process if the following conditional expectations satisfy the
relation}%
\begin{equation}
\mathbb{E}\left(  h(Z_{r})\left\vert \mathcal{F}_{s}^{+}\vee\mathcal{F}%
_{t}^{-}\right.  \right)  =\mathbb{E}\left(  h(Z_{r})\left\vert Z_{s}%
,Z_{t}\right.  \right)  \label{conditionalexpectations}%
\end{equation}
\textit{for every bounded Borel measurable function }$h:$\textit{ }%
$\overline{D}\mapsto\mathbb{R}$, \textit{and for all }$r,s,t$\textit{
satisfying }$r\in\left(  s,t\right)  \subset\left[  0,T\right]  $\textit{. In
(\ref{conditionalexpectations}), }$\mathcal{F}_{s}^{+}$ \textit{denotes the
}$\sigma$\textit{-algebra generated by the }$Z_{\tau}$\textit{'s for all
}$\tau\in\left[  0,s\right]  $, \textit{while }$\mathcal{F}_{t}^{-}$
\textit{is that generated by the }$Z_{\tau}$\textit{'s for all }$\tau
\in\left[  t,T\right]  $.

\bigskip

In the sequel we shall denote by $\mathcal{F}_{\tau\in\left[  0,T\right]
}^{+}$ the \textit{increasing} filtration generated by the $\mathcal{F}%
_{s}^{+}$'s, and by $\mathcal{F}_{\tau\in\left[  0,T\right]  }^{-}$ the
\textit{decreasing} filtration generated by the $\mathcal{F}_{t}^{-}$'s.

In order to construct such processes with relation to (\ref{parabolicproblem})
and (\ref{adjointproblem}) we first recast Schr\"{o}dinger's and Bernstein's
ideas to fit the theory developed in \cite{jamison}. Accordingly, the two main
ingredients we need are transition density functions for the processes
together with joint probability distributions for $Z_{0}$ and $Z_{T}$. As we
shall see, this requires the existence of \textit{classical, positive}
solutions to (\ref{parabolicproblem})\ and (\ref{adjointproblem}),
respectively, which in turn requires good smoothness properties of $k$, $l$,
$V$, $\varphi$ and $\psi$. In order to achieve this we assume that the
boundary $\partial D$ is of class $\mathcal{C}^{2+\alpha}$ for some $\alpha
\in\left(  0,1\right)  $, and then impose the following hypotheses where
$n(x)$ denotes the unit outer normal vector at $x\in\partial D$ (in all that
follows we write $c$ for all the irrelevant constants that occur in the
various estimates unless we specify these constants otherwise, and we refer to
\cite{solonnikov} for a definition of the above concepts and for various
properties of the spaces of H\"{o}lder continuous functions introduced here):

\bigskip

(K) The function $k:\overline{D}\times\left[  0,T\right]  \mapsto
\mathbb{R}^{d^{2}}$ is such that for every $i,j\in\left\{  1,...,d\right\}  $
we have $k_{i,j}=k_{j,i}\in\mathcal{C}^{\alpha,\frac{\alpha}{2}}(\overline
{D}\times\left[  0,T\right]  )$ and $\frac{\partial k_{i,j}}{\partial x_{l}%
}\in\mathcal{C}^{\alpha,\frac{\alpha}{2}}(\overline{D}\times\left[
0,T\right]  )$ for all $i,j,l$. Moreover, the uniform ellipticity condition%
\[
\left(  k(x,t)q,q\right)  _{\mathbb{R}^{d}}\geq\underline{k}\left\vert
q\right\vert ^{2}%
\]
with $\underline{k}>0$ holds for all $(x,t)\in\overline{D}\times\left[
0,T\right]  $ and all $q\in\mathbb{R}^{d}$, where $\left\vert .\right\vert $
denotes the Euclidean norm. Finally, the conormal vector-field $n_{k}%
(x,t):=k(x,t)n(x)$ is uniformly outward pointing, nowhere tangent to $\partial
D$ for every $t\in\left[  0,T\right]  $ and we have%
\[
(x,t)\mapsto\sum_{i=1}^{d}k_{i,j}(x,t)n_{i}(x)\in\mathcal{C}^{1+\alpha
,\frac{1+\alpha}{2}}(\partial D\times\left[  0,T\right]  )
\]
for each $j$.

\bigskip

(L) For the components of the vector-field $l:\overline{D}\times\left[
0,T\right]  \mapsto\mathbb{R}^{d}$ we have $l_{i},\frac{\partial l_{i}%
}{\partial x_{j}}\in\mathcal{C}^{\alpha,\frac{\alpha}{2}}(\overline{D}%
\times\left[  0,T\right]  )$ for all $i,j\in\left\{  1,...,d\right\}  $.

\bigskip

(V) The function $V:\overline{D}\times\left[  0,T\right]  \mapsto\mathbb{R}$
is such that $V\in\mathcal{C}^{\alpha,\frac{\alpha}{2}}(\overline{D}%
\times\left[  0,T\right]  )$.

\bigskip

Thus, the above functions are all jointly H\"{o}lder continuous in the
space-time variable $(x,t)$.

Finally, $\varphi$ and $\psi$ ought to be smooth enough as well and compatible
with the boundary conditions in (\ref{parabolicproblem}) and
(\ref{adjointproblem}):

\bigskip

(IF) We have $\varphi,\psi\in\mathcal{C}^{2+\alpha}(\overline{D})$ with
$\varphi$ satisfying the conormal boundary condition relative to $k$ at $t=0$,
and $\psi$ satisfying that condition at $t=T$.

\bigskip

It then follows from the classical theory of linear parabolic equations (see
for instance \cite{friedman}, or more specifically Chapter 4 in
\cite{ladyurasolo} and Theorem 1 in \cite{eideliva}), which is, of course, a
particular case of the variational approach reviewed in the appendix, that
there exist evolution systems $U_{A}(t,s)_{0\leq s\leq t\leq T}$ and
$U_{A}^{\ast}(t,s)_{0\leq s\leq t\leq T}$ in $L^{2}(D)$ given by%
\begin{equation}
U_{A}(t,s)f(x)=\left\{
\begin{array}
[c]{c}%
f(x)\text{ \ \ \ \ \ \ \ \ \ \ \ \ \ \ \ \ \ \ \ \ \ \ \ \ if }t=s,\\
\int_{D}dyg_{A}(x,t;y,s)f(y)\text{ \ \ if }t>s
\end{array}
\right.  \label{evolutionsystem1}%
\end{equation}
and%
\begin{equation}
U_{A}^{\ast}(t,s)f(x)=\left\{
\begin{array}
[c]{c}%
f(x)\text{ \ \ \ \ \ \ \ \ \ \ \ \ \ \ \ \ \ \ \ \ \ \ \ \ if }t=s,\\
\int_{D}dyg_{A}^{\ast}(x,s;y,t)f(y)\text{ \ \ if }t>s
\end{array}
\right.  \label{evolutionsystem2}%
\end{equation}
with $g_{A}$ and $g_{A}^{\ast}$ the parabolic Green functions associated with
(\ref{parabolicproblem}) and (\ref{adjointproblem}). These functions satisfy%
\begin{equation}
g_{A}^{\ast}(x,s;y,t)=g_{A}(y,t;x,s) \label{greenfunctions}%
\end{equation}
for all $s,t\in\left[  0,T\right]  $ with $t>s$, and furthermore the functions
defined by%
\begin{equation}
u_{\varphi}(x,t):=\int_{D}dyg_{A}(x,t;y,0)\varphi(y),\text{ \ \ }t\in\left(
0,T\right]  \label{forwardsolution}%
\end{equation}
and%
\begin{equation}
v_{\psi}(x,t):=\int_{D}dyg_{A}^{\ast}(x,t;y,T)\psi(y),\text{ \ \ }t\in\left[
0,T\right)  , \label{backwardsolution}%
\end{equation}
are indeed classical solutions to (\ref{parabolicproblem}) and
(\ref{adjointproblem}), respectively. More precisely we have the following
result for them:

\bigskip

\textbf{Proposition 1. }\textit{Assume that Hypotheses (K), (L), (V) and
(IF)\ hold. Then the following statements are valid:}

\textit{(a)}\textbf{ }\textit{We have} $u_{\varphi},v_{\psi}\in\mathcal{C}%
^{2+\alpha,1+\frac{\alpha}{2}}(\overline{D}\times\left[  0,T\right]  )$
\textit{with }$u_{\varphi}$\textit{ the unique classical solution to
(\ref{parabolicproblem}) and }$v_{\psi}$\textit{ the unique classical solution
to (\ref{adjointproblem}).}

\textit{(b) If }$\varphi>0$\textit{, }$\psi>0$\textit{ on }$\overline{D}%
$\textit{ we have }$u_{\varphi}>0$\textit{, }$v_{\psi}>0$\textit{ on
}$\overline{D}\times\left[  0,T\right]  $\textit{, respectively.}

\textit{(c) If }$\varphi>0$\textit{ on }$\overline{D}$\textit{ we have }%
$g_{A}>0$\textit{ for all }$x,y\in$\textit{ }$\overline{D}$ \textit{and all
}$s,t\in\left[  0,T\right]  $ \textit{with} $t>s$ \textit{and furthermore
}$g_{A}$\textit{ is jointly continuous in these variables. Moreover, this
function} \textit{is} \textit{twice continuously differentiable in }%
$x$\textit{, once continuously differentiable in }$t$ \textit{and satisfies }%
\begin{align}
\partial_{t}g_{A}(x,t;y,s)  &  =-A(t)g_{A}(x,t;y,s),\text{ \ \ }(x,t)\in
D\times\left(  s,T\right]  ,\label{differentialequation}\\
\frac{\partial g_{A}(x,t;y,s)}{\partial n_{k}(x,t)}  &  =0,\text{
\ \ }(x,t)\in\partial D\times\left(  s,T\right]  ,\nonumber
\end{align}
\textit{where the elliptic differential operator}%
\begin{equation}
A(t):=-\frac{1}{2}\operatorname{div}\left(  k(.,t)\nabla\right)  +\left(
l(.,t),\nabla\right)  _{\mathbb{R}^{d}}+V(.,t) \label{ellipticoperator}%
\end{equation}
\textit{corresponds to the right-hand side of (\ref{parabolicproblem}).}
\textit{Finally, the heat kernel estimate}%
\begin{equation}
g_{A}(x,t;y,s)\leq c\left(  t-s\right)  ^{-\frac{d}{2}}\exp\left[
-c\frac{\left\vert x-y\right\vert ^{2}}{t-s}\right]
\label{heatkernelestimate}%
\end{equation}
\textit{holds}.

\textit{(d) If }$\psi>0$\textit{ on }$\overline{D}$ \textit{we have }%
$g_{A}^{\ast}>0$ \textit{for all }$x,y\in$\textit{ }$\overline{D}$ \textit{and
all }$s,t\in\left[  0,T\right]  $ \textit{with} $t>s$ \textit{and furthermore
}$g_{A}^{\ast}$\textit{ is jointly continuous in these variables. Moreover,
this function is twice continuously differentiable in }$x$\textit{, once
continuously differentiable in }$s$ \textit{and satisfies}%
\begin{align}
-\partial_{s}g_{A}^{\ast}(x,s;y,t)  &  =-A^{\ast}(s)g_{A}^{\ast}%
(x,s;y,t),\text{ \ \ }(x,s)\in D\times\left[  0,t\right)
,\label{adjointdifferentialequation}\\
\frac{\partial g_{A}^{\ast}(x,s;y,t)}{\partial n_{k}(x,s)}  &  =0,\text{
\ \ }(x,s)\in\partial D\times\left[  0,t\right)  ,\nonumber
\end{align}
\textit{where}%
\begin{equation}
A^{\ast}(s):=-\frac{1}{2}\operatorname{div}\left(  k(.,s)\nabla\right)
-\operatorname{div}\left(  \cdot l(.,s)\right)  +V(.,s)
\label{adjointellipticoperator}%
\end{equation}
\textit{is the formal adjoint to (\ref{ellipticoperator}) corresponding to the
right-hand side of (\ref{adjointproblem}). Finally, the same heat kernel
estimate}%
\begin{equation}
g_{A}^{\ast}(x,s;y,t)\leq c\left(  t-s\right)  ^{-\frac{d}{2}}\exp\left[
-c\frac{\left\vert x-y\right\vert ^{2}}{t-s}\right]
\label{heatkernelestimate1}%
\end{equation}
\textit{ as in (c) holds.}

\bigskip

Hypotheses (K), (L), (V) and (IF) will be our standing hypotheses in the sequel.

Let us now consider the function%
\begin{equation}
P\left(  x,t;E,r;y,s\right)  :=\int_{E}dzp(x,t;z,r;y,s)
\label{transitionfunction1}%
\end{equation}
for every $E\in\mathcal{B(}\overline{D})$, where%
\begin{equation}
p(x,t;z,r;y,s):=\frac{g_{A}(x,t;z,r)g_{A}(z,r;y,s)}{g_{A}(x,t;y,s)}
\label{transitiondensity}%
\end{equation}
is well defined and positive for all $x,y,z\in\overline{D}$ and all $r,s,t$
satisfying $r\in\left(  s,t\right)  \subset\left[  0,T\right]  $ when the
first part of (c)\ in Proposition 1 holds. In addition, let $\mu$ be a
positive integrable function on $\overline{D}\times\overline{D}$ such that%
\begin{equation}
\mu(E\times F):=\int_{E\times F}dxdy\mu(x,y) \label{probabilitymeasure}%
\end{equation}
defines a probability measure on $\mathcal{B(}\overline{D})\times
\mathcal{B(}\overline{D})$. It is remarkable that the simultaneous knowledge
of (\ref{transitionfunction1}) and (\ref{probabilitymeasure}) is sufficient to
guarantee the existence of a Bernstein process associated with
(\ref{parabolicproblem}). Indeed we have the following result:

\bigskip

\textbf{Theorem 1. }\textit{Assume that the first part of (c) in Proposition 1
holds, and let }$\mu$ \textit{be given by (\ref{probabilitymeasure}).}
\textit{Then there exists a probability space} $\left(  \Omega,\mathcal{F}%
,\mathbb{P}_{\mu}\right)  $ \textit{and a }$\overline{D}$\textit{-valued
Bernstein process }$Z_{\tau\in\left[  0,T\right]  }$\textit{ on} $\left(
\Omega,\mathcal{F},\mathbb{P}_{\mu}\right)  $ \textit{such that}%
\begin{equation}
\mathbb{P}_{\mu}\left(  Z_{0}\in E_{0},Z_{T}\in E_{T}\right)  =\mu(E_{0}\times
E_{T}) \label{endpointprobability}%
\end{equation}
\textit{for all }$E_{0},E_{T}\in\mathcal{B(}\overline{D})$.\textit{
Furthermore we have}%
\[
\mathbb{P}_{\mu}\left(  Z_{r}\in E\left\vert Z_{s},Z_{t}\right.  \right)
=P\left(  Z_{t},t;E,r;Z_{s},s\right)
\]
\textit{for every }$E\in\mathcal{B(}\overline{D})$ \textit{and all }%
$r,s,t$\textit{ satisfying }$r\in\left(  s,t\right)  \subset\left[
0,T\right]  $\textit{, and moreover the finite-dimensional distributions are
given by}%
\begin{align}
&  \mathbb{P}_{\mu}\left(  Z_{0}\in E_{0},Z_{t_{1}}\in E_{1},...,Z_{t_{n}}\in
E_{n},Z_{T}\in E_{T}\right) \nonumber\\
&  =\int_{E_{0}\times E_{T}}dxdy\mu(x,y)\int_{E_{1}}dx_{1}...\int_{E_{n}%
}dx_{n}%
{\displaystyle\prod\limits_{i=1}^{n}}
p\left(  y,T;x_{i},t_{i};x_{i-1},t_{i-1}\right)  \label{pathprobabilaw}%
\end{align}
\textit{where }$x_{0}=x$\textit{,} \textit{for all }$E_{0},E_{1}%
,...,E_{n},E_{T}\in\mathcal{B}(\overline{D})$\textit{ and all }$t_{0}%
,...,t_{n}\in\left[  0,T\right)  $ \textit{satisfying }$t_{0}=0<t_{1}%
<...<t_{n}<T$. \textit{Finally, }$\mathbb{P}_{\mu}$\textit{ is the unique
probability measure with these properties.}

\bigskip

\textbf{Proof. }The mapping $(x,y)\mapsto P\left(  x,t;E,r;y,s\right)  $ is
evidently continuous on $\overline{D}\times\overline{D}$ for every
$E\in\mathcal{B}(\overline{D})$ and all $r,s,t$ satisfying $r\in\left(
s,t\right)  \subset\left[  0,T\right]  $. Moreover, the mapping $E\mapsto
P\left(  x,t;E,r;y,s\right)  $ defines a probability measure on $\mathcal{B}%
(\overline{D})$ for all $x,y\in\overline{D}$ and all of those $r,s,t$; this is
indeed a consequence of the composition law%
\begin{equation}
U_{A}(t,s)=U_{A}(t,r)U_{A}(r,s) \label{compositionlaw}%
\end{equation}
pertaining to the evolution system (\ref{evolutionsystem1}), which translates
as%
\begin{equation}
g_{A}(x,t;y,s)=\int_{D}dzg_{A}(x,t;z,r)g_{A}(z,r;y,s) \label{compositiongreen}%
\end{equation}
for the corresponding Green function. We also have the relation%
\begin{align}
&  \int_{E}dxp(x^{\prime},t^{\prime};x,t;y;s)P\left(  x,t;F,r;y,s\right)
\nonumber\\
&  =\int_{F}dzp(x^{\prime},t^{\prime};z,r;y;s)P\left(  x^{\prime},t^{\prime
};E,t;z,r\right)  \label{importantrelation}%
\end{align}
\newline for all $E,F\in\mathcal{B}(\overline{D})$, all $x^{\prime}%
,y\in\overline{D}$ and all $r,s,t,t^{\prime}$ satisfying $r\in\left(
s,t\right)  \subset\left(  s,t^{\prime}\right)  \subset\left[  0,T\right]  $.
In order to see this we remark that the relation%
\begin{align*}
&  p(x^{\prime},t^{\prime};x,t;y,s)p(x,t;z,r;y,s)\\
&  =p(x^{\prime},t^{\prime};z,r;y,s)p(x^{\prime},t^{\prime};x,t;z,r)
\end{align*}
holds as an immediate consequence of (\ref{transitiondensity}), so that
(\ref{importantrelation}) obtains by integrating both sides of the preceding
identity first with respect to $z$ over $F$, and then the resulting expression
with respect to $x$ over $E$. The statement of the theorem then follows from a
direct adaptation of the arguments in Section 2 of \cite{jamison}.
\ \ $\blacksquare$

\bigskip

\textsc{Remark. }The preceding result makes it clear that the probability of
having $Z_{r}\in E$ is indeed solely conditioned by the past information at
time $s$ and the future information at time $t$. Furthermore, owing to
(\ref{endpointprobability}) and (\ref{pathprobabilaw}) the probability measure
$\mu$ clearly plays the r\^{o}le of a joint endpoint distribution for
$Z_{0},Z_{T}$, and we ought to note that there is \textit{a priori} no reason
why $Z_{\tau\in\left[  0,T\right]  }$ should be a Markov process when $\mu$ is
arbitrary. However $Z_{\tau\in\left[  0,T\right]  }$ does become a reversible
Markov diffusion in a sense we shall define shortly, for a very special class
of endpoint distributions which we will identify. For this we first associate
a \textit{forward} Markov transition function with (\ref{adjointproblem}) in
the following way:

\bigskip

\textbf{Lemma 1.}\textit{ Assume that the part of (b) in Proposition 1
relative to }$\psi$ \textit{and }$v_{\psi}$\textit{ holds, together with the
first part of (d). Let us define the function}%
\begin{equation}
M^{\ast}\left(  x,s;E,t\right)  :=\int_{E}dym^{\ast}(x,s;y,t)
\label{forwardtransition}%
\end{equation}
\textit{for each }$E\in\mathcal{B}(\overline{D})$\textit{, every }%
$x\in\overline{D}$\textit{ and all }$s,t$\textit{ }$\in\left[  0,T\right]
$\textit{ with }$t>s$\textit{, where}%
\begin{equation}
m^{\ast}(x,s;y,t):=g_{A}^{\ast}(x,s;y,t)\frac{v_{\psi}(y,t)}{v_{\psi}(x,s)}
\label{markovdensity1}%
\end{equation}
\textit{with }$v_{\psi}$\textit{ given by (\ref{backwardsolution}). Then
(\ref{forwardtransition}) is the transition function of a forward Markov
process in }$\overline{D}$.

\bigskip

\textbf{Proof. }The mapping $x\mapsto M^{\ast}\left(  x,s;E,t\right)  $ is
positive and continuous on $\overline{D}$, and furthermore $E\mapsto M^{\ast
}\left(  x,s;E,t\right)  $ defines a probability measure on $\mathcal{B(}%
\overline{D})$. This last assertion follows from the relation%
\[
v_{\psi}(x,s)=\int_{D}dyg_{A}^{\ast}(x,s;y,t)v_{\psi}(y,t),
\]
which, in turn, is a simple consequence of the composition law%
\[
U_{A}^{\ast}(t,s)=U_{A}^{\ast}(r,s)U_{A}^{\ast}(t,r)
\]
for the evolution system (\ref{evolutionsystem2}), which gives%
\begin{equation}
g_{A}^{\ast}(x,s;y,t)=\int_{D}dzg_{A}^{\ast}(x,s;z,r)g_{A}^{\ast}(z,r;y,t)
\label{adjointcompositiongreen}%
\end{equation}
for the corresponding Green function. But then the Chapman-Kolmogorov relation%
\[
M^{\ast}\left(  x,s;E,t\right)  =\int_{D}dym^{\ast}(x,s;y,r)M^{\ast}\left(
y,r;E,t\right)
\]
holds for (\ref{forwardtransition}), since (\ref{markovdensity1}) and
(\ref{adjointcompositiongreen}) imply that%
\[
m^{\ast}(x,s;y,t)=\int_{D}dzm^{\ast}(x,s;z,r)m^{\ast}(z,r;y,t)
\]
for all $r,s,t$ such that $r\in\left(  s,t\right)  \subset\left[  0,T\right]
$. \ \ $\blacksquare$

\bigskip

In a completely symmetric way, we can also associate a \textit{backward}
Markov transition function with (\ref{parabolicproblem}). Indeed we have the
following result, whose proof is entirely similar to that of the preceding
lemma and thereby omitted:

\bigskip

\textbf{Lemma 2. }\textit{Assume that the part of (b) in Proposition 1
relative to }$\varphi$ \textit{and }$u_{\varphi}$\textit{ holds, together with
the first part of (c). Let us define the function}%
\begin{equation}
M\left(  x,t;E,s\right)  :=\int_{E}dym(x,t;y,s) \label{backwardtransition}%
\end{equation}
\textit{for each }$E\in\mathcal{B}(\overline{D})$\textit{, every }%
$x\in\overline{D}$\textit{ and all }$s,t$\textit{ }$\in\left[  0,T\right]
$\textit{ with }$t>s$\textit{, where}%
\begin{equation}
m(x,t;y,s):=g_{A}(x,t;y,s)\frac{u_{\varphi}(y,s)}{u_{\varphi}(x,t)}
\label{markovdensity2}%
\end{equation}
\textit{with }$u_{\varphi}$\textit{ given by (\ref{forwardsolution}). Then
(\ref{backwardtransition}) is the transition function of a backward Markov
process in }$\overline{D}$.

\bigskip

The remarkable fact is that when $Z_{\tau\in\left[  0,T\right]  \text{ }}$is
reversible in the sense of the definition below, it becomes a realization of
the Markov processes we are alluding to in the preceding two lemmas. The
precise notion we need is the following, where we recall that $\varphi>0$,
$\psi>0$:

\bigskip

\textsc{Definition 2.} \textit{We say the Bernstein process }$Z_{\tau
\in\left[  0,T\right]  }$\textit{ of Theorem 1 is reversible if the density of
the joint probability measure (\ref{probabilitymeasure}) is of the form}%
\begin{equation}
\mu(x,y)=\varphi(x)g_{A}(y,T;x,0)\psi(y) \label{specialclass}%
\end{equation}
\textit{where}%
\begin{equation}
\int_{D\times D}dxdy\varphi(x)g_{A}(y,T;x,0)\psi(y)=1.
\label{normalizationcondition}%
\end{equation}

\bigskip

For the corresponding initial and final marginal distributions we then have%
\[
\mu_{0}(E):=\mu(E\times D)=\int_{E}dx\varphi(x)v_{\psi}(x,0)
\]
and%
\[
\mu_{T}(F):=\mu(D\times F)=\int_{F}dy\psi(y)u_{\varphi}(y,T)
\]
as a consequence of (\ref{backwardsolution}) and (\ref{forwardsolution}),
respectively, with $\mu_{0}(D)=\mu_{T}(D)=1$. Our definition is motivated by
the following result, where we use the shorthand notation%
\begin{equation}
\mu_{0}(x)=\varphi(x)v_{\psi}(x,0) \label{initialmarginaldensity}%
\end{equation}
and%
\begin{equation}
\mu_{T}(y)=\psi(y)u_{\varphi}(y,T) \label{finalmarginaldensity}%
\end{equation}
for the marginal densities:

\bigskip

\textbf{Theorem 2. }\textit{Assume that (b) and the first parts of (c) and (d)
in Proposition 1 hold. Assume furthermore that }$\mu$ \textit{is given by
(\ref{probabilitymeasure}) and (\ref{specialclass}), and let }$Z_{\tau
\in\left[  0,T\right]  }$\textit{ be the corresponding Bernstein process
associated with (\ref{parabolicproblem}) in the sense of Theorem 1. Then the
following statements hold for its finite-dimensional distributions:}

\textit{(a) We have}%
\begin{align}
&  \mathbb{P}_{\mu}\left(  Z_{0}\in E_{0},Z_{t_{1}}\in E_{1},...,Z_{t_{n}}\in
E_{n}\right) \nonumber\\
&  =\int_{E_{0}}dx\mu_{0}(x)\int_{E_{1}}dx_{1}...\int_{E_{n}}dx_{n}%
{\displaystyle\prod\limits_{i=1}^{n}}
m^{\ast}\left(  x_{i-1},t_{i-1};x_{i},t_{i}\right)  \label{markovprobabilaw}%
\end{align}
\textit{where }$m^{\ast}$ \textit{is given by (\ref{markovdensity1}) and
}$x_{0}=x$\textit{,} \textit{for all }$E_{0},E_{1},...,E_{n}\in\mathcal{B}%
(\overline{D})$\textit{ and all }$t_{0},...,t_{n}\in\left[  0,T\right)  $
\textit{satisfying }$t_{0}=0<t_{1}<...<t_{n}<T$. \textit{Thus }$Z_{\tau
\in\left[  0,T\right]  }$\textit{ is a forward Markov process with transition
function }$m^{\ast}$ \textit{and initial distribution density }$\mu_{0}$.

\textit{(b) We have}%
\begin{align}
&  \mathbb{P}_{\mu}\left(  Z_{T}\in E_{T},Z_{t_{n}}\in E_{n},...,Z_{t_{1}}\in
E_{1}\right) \nonumber\\
&  =\int_{E_{T}}dy\mu_{T}(y)\int_{E_{1}}dx_{1}...\int_{E_{n}}dx_{n}%
{\displaystyle\prod\limits_{i=1}^{n}}
m\left(  x_{i+1},t_{i+1};x_{i},t_{i}\right)  \label{markovprobilawbis}%
\end{align}
\textit{where }$m$ \textit{is given by (\ref{markovdensity2}) and }$x_{n+1}%
=y$\textit{,} \textit{for all }$E_{T},E_{1},...,E_{n}\in\mathcal{B}%
(\overline{D})$\textit{ and all }$t_{1},...,t_{n+1}\in\left(  0,T\right]  $
\textit{satisfying }$T=t_{n+1}>t_{n}>...>t_{1}>0$. \textit{Thus }$Z_{\tau
\in\left[  0,T\right]  }$\textit{ is also a backward Markov process with
transition function }$m$ \textit{and final distribution density }$\mu_{T}$.

\textit{(c)} \textit{We have}%
\begin{equation}
\mathbb{P}_{\mu}\left(  Z_{t}\in E\right)  =\int_{E}dxu_{\varphi}(x,t)v_{\psi
}(x,t) \label{bernsteindensity}%
\end{equation}
\textit{for each }$E\in\mathcal{B}(\overline{D})$\textit{ and every }%
$t\in\left[  0,T\right]  $\textit{, where }$u_{\varphi}$\textit{ and }%
$v_{\psi}$\textit{ are given by (\ref{forwardsolution}) and
(\ref{backwardsolution}), respectively. }

\bigskip

\textbf{Proof. }From (\ref{pathprobabilaw}) with $E_{T}=D$ and
(\ref{transitiondensity}) we have%
\begin{align}
&  \mathbb{P}_{\mu}\left(  Z_{0}\in E_{0},Z_{t_{1}}\in E_{1},...,Z_{t_{n}}\in
E_{n}\right) \nonumber\\
&  =\int_{E_{0}\times D}dxdy\mu(x,y)\int_{E_{1}}dx_{1}...\int_{E_{n}}dx_{n}%
{\displaystyle\prod\limits_{i=1}^{n}}
\frac{g_{A}(y,T;x_{i},t_{i})g_{A}(x_{i},t_{i};x_{i-1},t_{i-1})}{g_{A}%
(y,T;x_{i-1},t_{i-1})}\nonumber\\
&  =\int_{E_{0}\times D}dxdy\varphi(x)\psi(y)\int_{E_{1}}dx_{1}...\int_{E_{n}%
}dx_{n}%
{\displaystyle\prod\limits_{i=1}^{n}}
g_{A}(x_{i},t_{i};x_{i-1},t_{i-1})\times g_{A}(y,T;x_{n},t_{n})\nonumber\\
&  =\int_{E_{0}}dx\varphi(x)\int_{E_{1}}dx_{1}...\int_{E_{n}}dx_{n}%
{\displaystyle\prod\limits_{i=1}^{n}}
g_{A}^{\ast}(x_{i-1},t_{i-1};x_{i},t_{i})\times v_{\psi}(x_{n},t_{n})
\label{markovprobabilaw1}%
\end{align}
after the use of (\ref{specialclass}), the successive cancellation of the
denominators in the above product and the use of (\ref{greenfunctions}) with
(\ref{backwardsolution}). Moreover, because of (\ref{markovdensity1}) we may
write%
\[
g_{A}^{\ast}(x_{i-1},t_{i-1};x_{i},t_{i})v_{\psi}(x_{i},t_{i})=v_{\psi
}(x_{i-1},t_{i-1})m^{\ast}\left(  x_{i-1},t_{i-1};x_{i},t_{i}\right)
\]
for every $i\in\left\{  1,...,n\right\}  $, so that the repeated application
of this relation in the product on the very right-hand side of
(\ref{markovprobabilaw1}) leads to%
\begin{align*}
&  \mathbb{P}_{\mu}\left(  Z_{0}\in E_{0},Z_{t_{1}}\in E_{1},...,Z_{t_{n}}\in
E_{n}\right) \\
&  =\int_{E_{0}}dx\varphi(x)v_{\psi}(x,0)\int_{E_{1}}dx_{1}...\int_{E_{n}%
}dx_{n}%
{\displaystyle\prod\limits_{i=1}^{n}}
m^{\ast}\left(  x_{i-1},t_{i-1};x_{i},t_{i}\right)  ,
\end{align*}
which is the desired result. The proof of (\ref{markovprobilawbis}) follows
from entirely similar arguments based on (\ref{forwardsolution}) and
(\ref{markovdensity2}), with $E_{0}=D$. Finally, (\ref{bernsteindensity}) is a
straightforward consequence of (\ref{markovprobabilaw}) with $E_{0}=D$, or of
(\ref{markovprobilawbis}) with $E_{T}=D$, in both cases with $n=1$.
\ \ $\blacksquare$

\bigskip

\textsc{Remarks.} (1) It follows from (\ref{adjointproblem}),
(\ref{adjointdifferentialequation}), (\ref{adjointellipticoperator}) and some
lengthy calculations that the transition density (\ref{markovdensity1})
satisfies the parabolic partial differential equation%
\begin{align}
&  -\partial_{s}m^{\ast}(x,s;y,t)\nonumber\\
&  =\frac{1}{2}\operatorname{div}_{x}\left(  k(x,s)\nabla_{x}m^{\ast
}(x,s;y,t)\right)  +\left(  l(x,s),\nabla_{x}m^{\ast}(x,s;y,t)\right)
_{\mathbb{R}^{d}}\nonumber\\
&  +\left(  k(x,s)\nabla_{x}\ln v_{\psi}(x,s),\nabla_{x}m^{\ast}%
(x,s;y,t)\right)  _{\mathbb{R}^{d}} \label{parabolicdensity}%
\end{align}
relative to the variables $(x,s)\in D\times\left[  0,t\right)  $ of the past,
along with the boundary condition%
\begin{equation}
\frac{\partial m^{\ast}(x,s;y,t)}{\partial n_{k}(x,s)}=0,\text{ \ \ }%
(x,s)\in\partial D\times\left[  0,t\right)  . \label{boundarycondition1}%
\end{equation}
In a similar way we infer from (\ref{parabolicproblem}),
(\ref{differentialequation}) and (\ref{ellipticoperator}) that the transition
density (\ref{markovdensity2}) satisfies the equation%
\begin{align}
&  \partial_{t}m(x,t;y,s)\nonumber\\
&  =\frac{1}{2}\operatorname{div}_{x}\left(  k(x,t)\nabla_{x}%
m(x,t;y,s)\right)  -\left(  l(x,t),\nabla_{x}m(x,t;y,s)\right)  _{\mathbb{R}%
^{d}}\nonumber\\
&  +\left(  k(x,t)\nabla_{x}\ln u_{\varphi}(x,t),\nabla_{x}m(x,t;y,s)\right)
_{\mathbb{R}^{d}} \label{parabolicdensitybis}%
\end{align}
with respect to the variables $(x,t)\in D\times\left(  s,T\right]  $ of the
future, together with the boundary condition%
\begin{equation}
\frac{\partial m(x,t;y,s)}{\partial n_{k}(x,t)}=0,\text{ \ \ }(x,t)\in\partial
D\times\left(  s,T\right]  . \label{boundarycondition2}%
\end{equation}
The preceding relations suggest that we may think of the reversible Markov
process of Theorem 2 as a process wandering in $\overline{D}$ which becomes
reflected in the conormal direction whenever it hits the boundary $\partial
D$. This way of looking at $Z_{\tau\in\left[  0,T\right]  }$ is reminiscent of
the definition of the standard reflected Brownian motion given at the very
beginning of \cite{hsu}, and we can indeed prove the reflection property we
just alluded to in the first example of Section 2. Further below we also
explain the appearance of the somewhat exotic logarithmic terms in
(\ref{parabolicdensity}) and (\ref{parabolicdensitybis}), whose structure is
intimately tied up with the specific form of (\ref{markovdensity1}) and
(\ref{markovdensity2}).

(2) The considerations that lead to the statement of Theorem 2 show that the
marginal densities (\ref{initialmarginaldensity}) and
(\ref{finalmarginaldensity}) are entirely determined by $\varphi$ and $\psi$.
We could also have adopted the inverse point of view, namely, that of
prescribing continuous $\mu_{0}>0$ and $\mu_{T}>0$ satisfying the
normalization conditions%
\[
\int_{D}dx\mu_{0}(x)=\int_{D}dx\mu_{T}(x)=1,
\]
and then considered the relations%
\begin{align}
\varphi(x)\int_{D}dzg_{A}(z,T;x,0)\psi(z)  &  =\mu_{0}(x),\nonumber\\
\psi(y)\int_{D}dzg_{A}(y,T;z,0)\varphi(z)  &  =\mu_{T}(y)
\label{integralequations}%
\end{align}
as a nonlinear inhomogeneous system of integral equations in the two unknowns
$\varphi$ and $\psi$. However, whereas it is true that such $\mu_{0}$ and
$\mu_{T}$ imply the existence of a unique solution to (\ref{integralequations}%
) consisting of continuous and positive functions $\varphi$ and $\psi$ as a
consequence of the main theorem in \cite{beurling}, the theory developed in
that article guarantees neither their regularity nor their the H\"{o}lder
continuity, which will be so crucial to our considerations. It is in fact an
interesting open problem whether the main result of \cite{beurling} can be
extended to cover such situations.

(3) Theorem 2 clearly illustrates a kind of reversibility of $Z_{\tau
\in\left[  0,T\right]  }$ in that this process can run back and forth within
$\overline{D}$, which \textit{a posteriori }justifies the terminology of
Definition 2. In particular, the probability density%
\begin{equation}
\rho_{\mu}(x,t):=u_{\varphi}(x,t)v_{\psi}(x,t) \label{probabilitydensity}%
\end{equation}
in (\ref{bernsteindensity}) is expressed as the product of solutions to
(\ref{parabolicproblem}) and (\ref{adjointproblem}), which indeed brings in
the two time directions in an explicit way. Therefore, at this stage it is
natural to start exploring the possible connections that might exist between
the preceding considerations and the notion of reversibility put forward in
\cite{kolmogorov}. On the one hand, it follows from (\ref{greenfunctions}),
(\ref{markovdensity1}) and (\ref{markovdensity2}) that the identity%

\begin{equation}
m(y,t;x,s)\rho_{\mu}(y,t)=\rho_{\mu}(x,s)m^{\ast}(x,s;y,t)
\label{reciprocalidentity}%
\end{equation}
holds for all $(x,y)\in$ $\overline{D}\times\overline{D}$\textit{ }and all
$s,t$ $\in\left[  0,T\right]  $ with $t>s$, which plays the r\^{o}le of (8) in
\cite{kolmogorov}. From (\ref{reciprocalidentity}) we then immediately infer
that%
\[
\rho_{\mu}(y,t)=\int_{D}dx\rho_{\mu}(x,s)m^{\ast}(x,s;y,t)
\]
and%
\[
\rho_{\mu}(x,s)=\int_{D}dym(y,t;x,s)\rho_{\mu}(y,t),
\]
which generalize (7) in \cite{kolmogorov}. On the other hand, however, let us
assume momentarily that $Z_{\tau\in\left[  0,T\right]  }$ is also reversible
in the sense of \cite{kolmogorov}, which means that%
\begin{equation}
m(y,t;x,s)=m^{\ast}(y,s;x,t) \label{kolmoreversible}%
\end{equation}
according to (9) of that article. From the preceding relation it then follows
at once from (\ref{greenfunctions}) that the equality%

\[
g_{A}(y,t;x,s)\frac{u_{\varphi}(x,s)}{u_{\varphi}(y,t)}=g_{A}(x,t;y,s)\frac
{v_{\psi}(x,t)}{v_{\psi}(y,s)}%
\]
is valid for all $(x,y)\in$ $\overline{D}\times\overline{D}$\textit{ }and all
$s,t$ $\in\left[  0,T\right]  $ with $t>s$, which implies that%
\[
\rho_{\mu}(x,t)=\rho_{\mu}(x,s)
\]
for each $x\in\overline{D}$ and every $t\geq s$ by choosing $y=x$. Therefore
$\rho_{\mu}$ must be independent of time, a very particular situation indeed
which is almost never realized in our context with the exception of a few
cases. Thus, whereas $Z_{\tau\in\left[  0,T\right]  }$ is a reversible Markov
process in the sense of Theorem 2 thanks to the very specific form of the
joint measures (\ref{specialclass}), it is in general not reversible according
to \cite{kolmogorov}. We will dwell more on this further in this section when
we have additional information about $Z_{\tau\in\left[  0,T\right]  }$.

(4) While a separable version of the process $Z_{\tau\in\left[  0,T\right]  }$
always exists under the above hypotheses, we may also assume that $Z_{\tau
\in\left[  0,T\right]  }$ is continuous. Indeed for all $s,t\in\left[
0,T\right]  $ with $t>s$ and $\gamma\in\left(  0,+\infty\right)  $ we have%
\[
\mathbb{E}_{\mu}\left\vert Z_{t}-Z_{s}\right\vert ^{\gamma}=\int_{D\times
D}dxdy\mu_{s}(x)m^{\ast}(x,s;y,t)\left\vert y-x\right\vert ^{\gamma}%
\]
where $\mu_{s}$ is the distribution density at time $s$, namely,%
\[
\mu_{s}(x)=\int_{D}dz\mu_{0}(z)m^{\ast}(z,0;x,s),
\]
with $\mathbb{E}_{\mu}$ the expectation functional with respect to
$\mathbb{P}_{\mu}$. Consequently, since $(x,s,y,t)\mapsto$ $v_{\psi}%
^{-1}(x,s)v_{\psi}(y,t)$ is uniformly bounded on $\overline{D}\times\left[
0,T\right]  \times\overline{D}\times\left[  0,T\right]  $ and since%
\[
\int_{D}dx\mu_{s}(x)=1,
\]
we get the estimate%
\begin{align}
&  \mathbb{E}_{\mu}\left\vert Z_{t}-Z_{s}\right\vert ^{\gamma}\nonumber\\
&  \leq c\int_{D}dx\mu_{s}(x)\int_{D}dyg_{A}^{\ast}(x,s;y,t)\left\vert
y-x\right\vert ^{\gamma}\nonumber\\
&  \leq c\left(  t-s\right)  ^{-\frac{d}{2}}\int_{D}dx\mu_{s}(x)\int_{D}%
dy\exp\left[  -c\frac{\left\vert y-x\right\vert ^{2}}{t-s}\right]  \left\vert
y-x\right\vert ^{\gamma}\nonumber\\
&  \leq c\left(  t-s\right)  ^{\frac{\gamma}{2}}\int_{\mathbb{R}^{d}}%
dy\exp\left[  -c\left\vert y\right\vert ^{2}\right]  \left\vert y\right\vert
^{\gamma}=c\left(  t-s\right)  ^{\frac{\gamma}{2}} \label{continuousversion}%
\end{align}
according to (\ref{heatkernelestimate1}), translation invariance on
$\mathbb{R}^{d}$ and an elementary change of variables. Therefore, fixing
$\gamma\in\left(  2,+\infty\right)  $ we obtain%
\[
\mathbb{E}_{\mu}\left\vert Z_{t}-Z_{s}\right\vert ^{\gamma}\leq c\left(
t-s\right)  ^{1+\delta}%
\]
with $\delta=\frac{\gamma-2}{2}$, so that the assertion follows from
Kolmogorov's continuity conditions. In the sequel we shall thereby always
assume that $Z_{\tau\in\left[  0,T\right]  }$ is separable and continuous.

\bigskip

While Theorem 2 shows that the process $Z_{\tau\in\left[  0,T\right]  }$ can
run forward and backward in a Markovian manner within $\overline{D}$ for the
very specific class of endpoint distributions (\ref{specialclass}), we now
proceed to investigate its dynamical properties more in detail. Our first step
in this direction is to show that $Z_{\tau\in\left[  0,T\right]  }$ is a
reversible Markov diffusion in $\overline{D}$.

\bigskip

\textbf{Lemma 3.} \textit{Assume that the same hypotheses as in Lemma 1 are
valid. Then Lindeberg's condition}%
\begin{equation}
\lim_{t\rightarrow s_{+}}\left(  t-s\right)  ^{-1}\int_{\left\{  y\in
D:\left\vert x-y\right\vert >\varepsilon\right\}  }dym^{\ast}(x,s;y,t)=0
\label{lindeberg1}%
\end{equation}
\textit{holds} \textit{uniformly in }$x\in D$\textit{ for each }$s$
\textit{and} \textit{every sufficiently small }$\varepsilon>0$\textit{. In a
similar way, if the same hypotheses as in Lemma 2 hold we have}%
\begin{equation}
\lim_{s\rightarrow t_{-}}\left(  t-s\right)  ^{-1}\int_{\left\{  y\in
D:\left\vert x-y\right\vert >\varepsilon\right\}  }dym(x,t;y,s)=0
\label{lindeberg2}%
\end{equation}
\textit{uniformly in }$x\in D$\textit{ for each }$t$ \textit{and every
sufficiently small }$\varepsilon>0$\textit{.}

\bigskip

\textbf{Proof.} In order to get (\ref{lindeberg1}) we must prove that%
\[
\lim_{t\rightarrow s_{+}}(t-s)^{-1}\int_{\left\{  y\in D:\left\vert
x-y\right\vert >\varepsilon\right\}  }dyg_{A}^{\ast}(x,s;y,t)\frac{v_{\psi
}(y,t)}{v_{\psi}(x,s)}=0
\]
uniformly in $x\in D$ for each $s$ according to (\ref{markovdensity1}). The
key observation for this is that the positive function $(x,s,y,t)\mapsto$
$v_{\psi}^{-1}(x,s)v_{\psi}(y,t)$ is uniformly bounded as a consequence of its
smoothness on $\overline{D}\times\left[  0,T\right]  \times$ $\overline
{D}\times\left[  0,T\right]  $, so that estimate (\ref{heatkernelestimate1})
leads to%
\begin{align*}
0  &  \leq\left(  t-s\right)  ^{-1}\int_{\left\{  y\in D:\left\vert
x-y\right\vert >\varepsilon\right\}  }dyg_{A}^{\ast}(x,s;y,t)\frac{v_{\psi
}(y,t)}{v_{\psi}(x,s)}\\
&  \leq c\left(  t-s\right)  ^{-\frac{d+2}{2}}\int_{\left\{  y\in D:\left\vert
x-y\right\vert >\varepsilon\right\}  }dy\exp\left[  -c\frac{\left\vert
y-x\right\vert ^{2}}{t-s}\right] \\
&  \leq c\left(  t-s\right)  ^{-\frac{d+2}{2}}\exp\left[  -c\frac
{\varepsilon^{2}}{t-s}\right]  \rightarrow0
\end{align*}
as $t\rightarrow s_{+}$ uniformly in $x$ since $D$ is bounded, as desired. The
proof of (\ref{lindeberg2}) is evidently identical, and based on
(\ref{heatkernelestimate}). \ \ $\blacksquare$

\bigskip

The first half of the preceding lemma now allows us to prove the following result:

\bigskip

\textbf{Proposition 2. }\textit{Assume that the same hypotheses as in Lemma 1
are valid. Then} \textit{the following statements hold:}

\textit{(a) We have}%
\begin{align}
&  \lim_{t\rightarrow s_{+}}\left(  t-s\right)  ^{-1}\int_{\left\{  y\in
D:\left\vert x-y\right\vert \leq\varepsilon\right\}  }dym^{\ast}%
(x,s;y,t)\left(  y-x\right) \nonumber\\
&  =a^{\ast}(x,s)+k(x,s)\nabla_{x}\ln v_{\psi}(x,s) \label{drift}%
\end{align}
\textit{for each }$x\in D$ \textit{and every} $s$ \textit{independently of any
sufficiently small }$\varepsilon>0$\textit{, where the }$i^{th}$%
\textit{component of the vector-field} $a^{\ast}$ \textit{is}%
\begin{equation}
a_{i}^{\ast}(x,s)=\frac{1}{2}\operatorname{div}_{x}\left(  k_{i}(x,s)\right)
+l_{i}(x,s) \label{driftterm1}%
\end{equation}
\textit{for every }$i\in\left\{  1,...,d\right\}  $\textit{, with }%
$k_{i}(x,s)$\textit{ the }$i^{th}$\textit{row or column of the symmetric
matrix }$k(x,s)$.

\textit{(b) We have}%
\begin{equation}
\lim_{t\rightarrow s_{+}}\left(  t-s\right)  ^{-1}\int_{\left\{  y\in
D:\left\vert x-y\right\vert \leq\varepsilon\right\}  }dym^{\ast}%
(x,s;y,t)\left(  \left(  y-x\right)  \otimes\left(  y-x\right)  \right)
=k(x,s) \label{diffusion}%
\end{equation}
\textit{for each }$x\in D$ \textit{and every} $s$,\textit{ independently of
any sufficiently small }$\varepsilon>0$.

\bigskip

\textbf{Proof. }Owing to (\ref{markovdensity1}) and (\ref{lindeberg1}) it is
sufficient to prove that%
\begin{align}
&  \lim_{t\rightarrow s_{+}}\left(  t-s\right)  ^{-1}\int_{D}dyg_{A}^{\ast
}(x,s;y,t)\frac{v_{\psi}(y,t)}{v_{\psi}(x,s)}\left(  y-x\right) \nonumber\\
&  =a^{\ast}(x,s)+k(x,s)\nabla_{x}\ln v_{\psi}(x,s) \label{drift2}%
\end{align}
and%
\begin{equation}
\lim_{t\rightarrow s_{+}}\left(  t-s\right)  ^{-1}\int_{D}dyg_{A}^{\ast
}(x,s;y,t)\frac{v_{\psi}(y,t)}{v_{\psi}(x,s)}\left(  \left(  y-x\right)
\otimes\left(  y-x\right)  \right)  =k(x,s), \label{diffusion2}%
\end{equation}
respectively, since the functions $y\mapsto y-x$ and $y\mapsto\left(
y-x\right)  \otimes\left(  y-x\right)  $ are bounded on $D$. Thanks to the
differentiability properties of $v_{\psi}$ and the convexity of $D$ we first
write%
\begin{align}
\frac{v_{\psi}(y,t)}{v_{\psi}(x,s)}  &  =1+\frac{\left(  \nabla_{x}v_{\psi
}(x,t),y-x\right)  _{\mathbb{R}^{d}}}{v_{\psi}(x,s)}+\frac{\left(
\mathsf{H}_{v_{\psi}}(x^{\ast},t)\left(  y-x\right)  ,y-x\right)
_{\mathbb{R}^{d}}}{2v_{\psi}(x,s)}\nonumber\\
&  +(t-s)\frac{\partial_{s^{\ast}}v_{\psi}(x,s^{\ast})}{v_{\psi}(x,s)}
\label{expansion}%
\end{align}
as a consequence of the Taylor expansion for $v_{\psi}$, where $\mathsf{H}%
_{v_{\psi}}$ denotes the Hessian matrix relative to the spatial variable alone
with $x^{\ast}$ a point on the open line segment joining $x$ and $y$, and
where $s^{\ast}\in\left(  s,t\right)  $. The strategy of the proof then
amounts to estimating the various contributions to (\ref{drift2}) and
(\ref{diffusion2}) coming from (\ref{expansion}).

In order to explain the appearance of the vector-field $a^{\ast}$ on the
right-hand side of (\ref{drift2}) we begin by showing that%
\begin{align}
&  \lim_{t\rightarrow s_{+}}\left(  t-s\right)  ^{-1}\int_{D}dyg_{A}^{\ast
}(x,s;y,t)\left(  y_{i}-x_{i}\right) \nonumber\\
&  =\frac{1}{2}\operatorname{div}_{x}\left(  k_{i}(x,s)\right)  +l_{i}(x,s)
\label{contribution1}%
\end{align}
for every $i\in\left\{  1,...,d\right\}  $. Let us define the function
$f_{i}:D\mapsto\mathbb{R}$ by $f_{i}(y):=y_{i}-x_{i}$; since $f_{i}(x)=0$ we
have%
\begin{align*}
&  \lim_{t\rightarrow s_{+}}\left(  t-s\right)  ^{-1}\int_{D}dyg_{A}^{\ast
}(x,s;y,t)\left(  y_{i}-x_{i}\right) \\
&  =\lim_{t\rightarrow s_{+}}\left(  t-s\right)  ^{-1}\left(  \int_{D}%
dyg_{A}^{\ast}(x,s;y,t)f_{i}(y)-f_{i}(x)\right) \\
&  =\left(  -A^{\ast}(s)f_{i}\right)  (x)=\frac{1}{2}\operatorname{div}%
_{x}\left(  k_{i}(x,s)\right)  +l_{i}(x,s)
\end{align*}
according to (\ref{evolutionsystem2}) and an elementary calculation based on
(\ref{adjointellipticoperator}), so that (\ref{contribution1}) holds. A
similar calculation with the function $f_{i,j}:D\mapsto\mathbb{R}$ given by
$f_{i,j}(y):=\left(  y_{i}-x_{i}\right)  \left(  y_{j}-x_{j}\right)  $ leads
to
\begin{align}
&  \lim_{t\rightarrow s_{+}}\left(  t-s\right)  ^{-1}\int_{D}dyg_{A}^{\ast
}(x,s;y,t)\left(  y_{i}-x_{i}\right)  \left(  y_{j}-x_{j}\right) \nonumber\\
&  =\left(  -A^{\ast}(s)f_{i,j}\right)  (x)=k_{i,j}(x,s) \label{contribution2}%
\end{align}
for all $i,j\in\left\{  1,...,d\right\}  $ since $f_{i,j}(x)=0$, which allows
us to evaluate the contribution of the gradient term in (\ref{expansion}).
Indeed, if we substitute that term into the $i^{th}$component of the left-hand
side of (\ref{drift2}) and use (\ref{contribution2}) we obtain%
\begin{align}
&  \lim_{t\rightarrow s_{+}}\left(  t-s\right)  ^{-1}%
{\displaystyle\int\limits_{D}}
dyg_{A}^{\ast}(x,s;y,t)\frac{\left(  \nabla_{x}v_{\psi}(x,t),y-x\right)
_{\mathbb{R}^{d}}}{v_{\psi}(x,s)}\left(  y_{i}-x_{i}\right) \nonumber\\
&  =\sum_{j=1}^{d}\frac{\partial_{x_{j}}v_{\psi}(x,s)}{v_{\psi}(x,s)}%
\lim_{t\rightarrow s_{+}}\left(  t-s\right)  ^{-1}\int_{D}dyg_{A}^{\ast
}(x,s;y,t)\left(  y_{i}-x_{i}\right)  \left(  y_{j}-x_{j}\right) \nonumber\\
&  =\frac{\left(  k(x,s)\nabla_{x}v_{\psi}(x,s)\right)  _{i}}{v_{\psi}(x,s)}
\label{contribution3}%
\end{align}
for every $i\in\left\{  1,...,d\right\}  $, which is the $i^{th}$component of
the second term on the right-hand side of (\ref{drift2}). Therefore, in order
to get (\ref{drift2}) it remains to prove that there are no contributions
coming from the third and fourth terms on the right-hand side of
(\ref{expansion}).

Regarding the third term we first observe that%
\[
\lim_{t\rightarrow s_{+}}\left(  t-s\right)  ^{-1}\int_{D}dyg_{A}^{\ast
}(x,s;y,t)\left\vert y-x\right\vert ^{3}=0
\]
uniformly in $x\in D$ for every $s$, since from (\ref{heatkernelestimate1}) we
obtain%
\begin{align*}
0  &  \leq\left(  t-s\right)  ^{-1}\int_{D}dyg_{A}^{\ast}(x,s;y,t)\left\vert
y-x\right\vert ^{3}\\
&  \leq c\left(  t-s\right)  ^{-\frac{d+2}{2}}\int_{D}dy\exp\left[
-c\frac{\left\vert y-x\right\vert ^{2}}{t-s}\right]  \left\vert y-x\right\vert
^{3}%
\end{align*}%
\begin{align}
&  \leq c\left(  t-s\right)  ^{-\frac{d+2}{2}}\int_{\mathbb{R}^{d}}%
dy\exp\left[  -c\frac{\left\vert y\right\vert ^{2}}{t-s}\right]  \left\vert
y\right\vert ^{3}\nonumber\\
&  =c(t-s)^{\frac{1}{2}}\int_{\mathbb{R}^{d}}dy\exp\left[  -c\left\vert
y\right\vert ^{2}\right]  \left\vert y\right\vert ^{3}\rightarrow0
\label{vanishingcontribution}%
\end{align}
as $t\rightarrow s_{+}$, again by translation invariance and the same change
of variables as in (\ref{continuousversion}). Consequently we have \textit{a
fortiori} the estimate%
\begin{align*}
0  &  \leq\left(  t-s\right)  ^{-1}\int_{D}dyg_{A}^{\ast}(x,s;y,t)\left\vert
\left(  \mathsf{H}_{v_{\psi}}(x^{\ast},t)\left(  y-x\right)  ,y-x\right)
_{\mathbb{R}^{d}}\left(  y_{i}-x_{i}\right)  \right\vert \\
&  \leq c\left(  t-s\right)  ^{-1}\int_{D}dyg_{A}^{\ast}(x,s;y,t)\left\vert
y-x\right\vert ^{3}\rightarrow0
\end{align*}
as $t\rightarrow s_{+}$ since the matrix-norm of $\mathsf{H}_{v_{\psi}%
}(x^{\ast},t)$ is uniformly bounded on the compact cylinder $\overline
{D}\times\left[  0,T\right]  $, from which we infer that%
\[
\lim_{t\rightarrow s_{+}}\left(  t-s\right)  ^{-1}%
{\displaystyle\int\limits_{D}}
dyg_{A}^{\ast}(x,s;y,t)\frac{\left(  \mathsf{H}_{v_{\psi}}(x^{\ast},t)\left(
y-x\right)  ,y-x\right)  _{\mathbb{R}^{d}}}{2v_{\psi}(x,s)}(y_{i}-x_{i})=0
\]
for each $s$ and every $i\in\left\{  1,...,d\right\}  $, as desired.

Finally, there is no contribution from the fourth term either since from its
substitution into the $i^{th}$component of the left-hand side of
(\ref{drift2}), the cancellation of the time increments and
(\ref{evolutionsystem2}) we get%
\begin{align*}
&  \frac{\partial_{t}v_{\psi}(x,s)}{v_{\psi}(x,s)}\lim_{t\rightarrow s_{+}%
}\int_{D}dyg_{A}^{\ast}(x,s;y,t)f_{i}(y)\\
&  =\frac{\partial_{t}v_{\psi}(x,s)}{v_{\psi}(x,s)}f_{i}(x)=0.
\end{align*}

The proof of (\ref{diffusion2}) is entirely similar, the unique non-vanishing
contribution being determined by the constant term on the right-hand side of
(\ref{expansion}) via (\ref{contribution2}). \ \ $\blacksquare$

\bigskip

The second half of Lemma 3 allows us to obtain a similar result for the
backward Markov process of Lemma 2. We omit the proof, which is based on
(\ref{differentialequation})-(\ref{heatkernelestimate}) of Proposition 1,
Lindeberg's condition (\ref{lindeberg2}), and is thereby identical to that of
the preceding proposition.

\bigskip

\textbf{Proposition 3. }\textit{Assume that the same hypotheses as in Lemma 2
are valid. Then} \textit{the following statements hold:}

\textit{(a) We have}%
\begin{align}
&  \lim_{s\rightarrow t_{-}}\left(  t-s\right)  ^{-1}\int_{\left\{  y\in
D:\left\vert x-y\right\vert \leq\varepsilon\right\}  }dym(x,t;y,s)\left(
x-y\right) \nonumber\\
&  =a(x,t)-k(x,t)\nabla_{x}\ln u_{\varphi}(x,t) \label{drift3}%
\end{align}
\textit{for each }$x\in D$ \textit{and every} $t$ \textit{independently of any
sufficiently small }$\varepsilon>0$\textit{, where the }$i^{th}$%
\textit{component of the vector-field} $a$ \textit{is}%
\begin{equation}
a_{i}(x,t)=-\frac{1}{2}\operatorname{div}_{x}\left(  k_{i}(x,t)\right)
+l_{i}(x,t) \label{driftterm2}%
\end{equation}
\textit{for every }$i\in\left\{  1,...,d\right\}  $\textit{, with }%
$k_{i}(x,t)$\textit{ as in Proposition 2}.

\textit{(b) We have}%
\begin{equation}
\lim_{s\rightarrow t_{-}}\left(  t-s\right)  ^{-1}\int_{\left\{  y\in
D:\left\vert x-y\right\vert \leq\varepsilon\right\}  }dym(x,t;y,s)\left(
\left(  x-y\right)  \otimes\left(  x-y\right)  \right)  =k(x,t)
\label{diffusion3}%
\end{equation}
\textit{for each }$x\in D$ \textit{and every} $t$,\textit{ independently of
any sufficiently small }$\varepsilon>0$.

\bigskip

Thus, both Propositions 2 and 3 show that the process $Z_{\tau\in\left[
0,T\right]  }$ of Theorem 2 is indeed a reversible Markov diffusion whose
coefficients are determined by (\ref{drift}), (\ref{diffusion}),
(\ref{drift3}) and (\ref{diffusion3}), respectively. In the sequel we shall
denote by
\begin{equation}
b^{\ast}(x,t):=a^{\ast}(x,t)+k(x,t)\nabla_{x}\ln v_{\psi}(x,t) \label{drift4}%
\end{equation}
and%
\begin{equation}
b(x,t):=a(x,t)-k(x,t)\nabla_{x}\ln u_{\varphi}(x,t) \label{drift5}%
\end{equation}
the drift terms (\ref{drift})\ and (\ref{drift3}), respectively.

\bigskip

\textsc{Remark.} We now explain why the specific form of (\ref{drift4}) and
(\ref{drift5}) suggests that our notion of reversibility and the results
obtained thus far may be interpreted as corresponding to a generalization of
the notion of reversibility defined in \cite{kolmogorov} which has its origins
in the last section of \cite{schroedinger}, and why our notion is suitable for
the description of non-stationary Markov processes such as $Z_{\tau\in\left[
0,T\right]  }$. In order to see this it is best to refer to the reformulation
of the main result of \cite{kolmogorov} as encoded in Relation (2.9) of
\cite{dobsukfritz}: given a time-homogeneous forward Markov diffusion
wandering in $\mathbb{R}^{d}$, a probability measure which is absolutely
continuous with respect to Lebesgue measure and defined by a smooth positive
density $\rho$ is reversible with respect to that diffusion if, and only if,
the generator of the diffusion is of the form%
\begin{align}
&  Lf(x)\nonumber\\
&  =\frac{1}{2\rho(x)}\operatorname{div}_{x}\left(  \rho(x)k(x)\nabla
_{x}f(x)\right) \nonumber\\
&  =\frac{1}{2}\operatorname{div}_{x}\left(  k(x)\nabla_{x}f(x)\right)
+\frac{1}{2}\left(  k(x)\nabla_{x}\ln\rho(x),\nabla_{x}f(x)\right)
_{\mathbb{R}^{d}} \label{dsfgenerator}%
\end{align}
for a suitable class of $f$'s, where $k$ is the associated symmetric, positive
definite and time-independent diffusion matrix. But this leads at once to%
\begin{equation}
b_{DSF}(x):=a_{DSF}(x)+\frac{1}{2}k(x)\nabla_{x}\ln\rho(x) \label{dsfdrift}%
\end{equation}
for the corresponding drift, with the $i^{th}$component of the vector-field
$a_{DSF}$ given by
\begin{equation}
a_{DSF,i}(x):=\frac{1}{2}\operatorname{div}_{x}\left(  k_{i}(x)\right)
\label{dsfvectorfield}%
\end{equation}
for every $i\in\left\{  1,...,d\right\}  $ where $k_{i}(x)$ is the $i^{th}$
row or column of $k(x)$. Indeed this follows immediately from
(\ref{dsfgenerator}) and the relation%
\begin{align*}
&  \frac{1}{2}\operatorname{div}_{x}\left(  k(x)\nabla_{x}f(x)\right) \\
&  =\frac{1}{2}\left(  k(x)\nabla_{x},\nabla_{x}f(x)\right)  _{\mathbb{R}^{d}%
}\text{ }+\left(  a_{DSF}(x),\nabla_{x}f(x)\right)  _{\mathbb{R}^{d}}\text{,}%
\end{align*}
so that in the end reversibility in the sense of \cite{dobsukfritz} or
\cite{kolmogorov} is equivalent to a drift of the form (\ref{dsfdrift}). But
it is then plain that (\ref{drift4}) and (\ref{drift5}) have the very same
structure as (\ref{dsfdrift}), and furthermore that%
\begin{align}
&  \frac{1}{2}\left(  b^{\ast}(x,t)-b(x,t)\right) \nonumber\\
&  =\hat{a}(x,t)+\frac{1}{2}k(x,t)\nabla_{x}\ln\rho_{\mu}(x,t)
\label{combinationofdrifts}%
\end{align}
where $\rho_{\mu}(x,t)$ is given by (\ref{probabilitydensity}) and%
\begin{equation}
\hat{a}_{i}(x,t)=\frac{1}{2}\operatorname{div}_{x}\left(  k_{i}(x,t)\right)
\label{vzvectorfield}%
\end{equation}
for every $i$, with (\ref{combinationofdrifts}) and (\ref{vzvectorfield})
formally identical to (\ref{dsfdrift}) and (\ref{dsfvectorfield}),
respectively. It is, therefore, deemed appropriate to say that the
reversibility of $Z_{\tau\in\left[  0,T\right]  }$ as illustrated by the
statement of Theorem 2 and due to the very specific form of the joint
probability measures (\ref{specialclass}) is a natural generalization of the
notion defined in \cite{dobsukfritz} and \cite{kolmogorov}. For another
discussion of reversibility in a more geometric context we refer the reader to
Section 4 in Chapter 5 of \cite{ikedawatanabe}.

\bigskip

We now proceed by proving that there exist two vector-valued Wiener processes
$W_{\tau\in\left[  0,T\right]  }^{\ast}$ and $W_{\tau\in\left[  0,T\right]  }%
$, indeed one for each one of the filtrations $\mathcal{F}_{\tau\in\left[
0,T\right]  }^{+}$ and $\mathcal{F}_{\tau\in\left[  0,T\right]  }^{-}$ we
defined at the very beginning of this section, which will eventually allow us
to consider $Z_{\tau\in\left[  0,T\right]  }$ as a forward and backward
It\^{o} diffusion whenever the drifts (\ref{drift4}) and (\ref{drift5}) do not
vanish identically and simultaneously. We begin with the following preparatory result:

\bigskip

\textbf{Lemma 4. }\textit{Assume that the hypotheses of Lemma 1 are valid and
let us define the process}%
\begin{equation}
Y_{t}^{\ast}:=Z_{t}-Z_{0}-\int_{0}^{t}d\tau b^{\ast}\left(  Z_{\tau}%
,\tau\right)  \label{continuousmartingale}%
\end{equation}
\textit{for every} $t\in\left[  0,T\right]  $\textit{, where} $Z_{\tau
\in\left[  0,T\right]  }$ \textit{is considered as the forward Markov
diffusion of Theorem 2. Then }$Y_{\tau\in\left[  0,T\right]  }^{\ast}$\textit{
is a continuous, square-integrable martingale with respect to} $\mathcal{F}%
_{\tau\in\left[  0,T\right]  }^{+}$\textit{. Under the hypotheses of Lemma 2 a
similar statement holds for the process}%
\begin{equation}
Y_{t}:=Z_{t}-Z_{T}+\int_{t}^{T}d\tau b\left(  Z_{\tau},\tau\right)
\label{continuousmartingale2}%
\end{equation}
\textit{with respect to} $\mathcal{F}_{\tau\in\left[  0,T\right]  }^{-}$,
\textit{with }$Z_{\tau\in\left[  0,T\right]  }$\textit{ considered as the
backward Markov diffusion of Theorem 2.}

\bigskip

\textbf{Proof. }We prove (\ref{continuousmartingale}) by first observing that%
\begin{equation}
\sup_{t\in\left[  0,T\right]  }\mathbb{E}_{\mu}\left\vert Z_{t}\right\vert
^{2}<+\infty\label{squareintegrability}%
\end{equation}
where%
\[
\mathbb{E}_{\mu}\left\vert Z_{t}\right\vert ^{2}=\int_{D\times D}dxdy\mu
_{0}(x)m^{\ast}(x,0;y,t)\left\vert y\right\vert ^{2}.
\]
Indeed we have%
\begin{align*}
&  \mathbb{E}_{\mu}\left\vert Z_{t}\right\vert ^{2}\\
&  \leq c\int_{D}dx\mu_{0}(x)\int_{D}dyg_{A}^{\ast}(x,0;y,t)\\
&  \leq ct^{-\frac{d}{2}}\int_{D}dx\mu_{0}(x)\int_{\mathbb{R}^{d}}%
dy\exp\left[  -c\frac{\left\vert y\right\vert ^{2}}{t}\right]  \leq c<+\infty
\end{align*}
uniformly in $t\in\left[  0,T\right]  $ according to
(\ref{heatkernelestimate1}), since $(x,y,t)\mapsto$ $v_{\psi}^{-1}%
(x,0)v_{\psi}(y,t)\left\vert y\right\vert ^{2}$ is uniformly bounded on
$\overline{D}\times\overline{D}\times\left[  0,T\right]  $. Therefore we have%
\begin{equation}
\sup_{t\in\left[  0,T\right]  }\mathbb{E}_{\mu}\left\vert Y_{t}^{\ast
}\right\vert ^{2}<+\infty\label{squareintegability2}%
\end{equation}
as a consequence of the uniform boundedness of $b^{\ast}(x,t)$ given by
(\ref{drift4}).

In order to prove the statement of the lemma, it is sufficient to show that
the scalar-valued process%
\begin{equation}
y_{t}^{\ast}:=\left(  Y_{t}^{\ast},q\right)  _{\mathbb{R}^{d}}=\left(
Z_{t}-Z_{0},q\right)  _{\mathbb{R}^{d}}-\int_{0}^{t}d\tau\left(  b^{\ast
}\left(  Z_{\tau},\tau\right)  ,q\right)  _{\mathbb{R}^{d}}
\label{scalarvaluedprocess}%
\end{equation}
is a continuous, square-integrable martingale for every $q\in$ $\mathbb{R}%
^{d}$. While the continuity is clear according to the remark following Theorem
2, the fact that%
\[
\sup_{t\in\left[  0,T\right]  }\mathbb{E}_{\mu}\left\vert y_{t}^{\ast
}\right\vert ^{2}<+\infty
\]
is an immediate consequence of (\ref{squareintegability2}) and
(\ref{scalarvaluedprocess}). Therefore, it remains to show that the equality%
\begin{equation}
\mathbb{E}_{\mu}\left(  y_{t}^{\ast}\left\vert \mathcal{F}_{s}^{+}\right.
\right)  =y_{s}^{\ast} \label{martingaleproperty}%
\end{equation}
holds $\mathbb{P}_{\mu}$-a.s. for all $t\geq s$, and for this we need only
prove that the right-hand derivative of (\ref{martingaleproperty}) with
respect to $t$ vanishes, namely, that%
\begin{equation}
\lim_{r\rightarrow t_{+}}\left(  r-t\right)  ^{-1}\mathbb{E}_{\mu}\left(
y_{r}^{\ast}-y_{t}^{\ast}\left\vert \mathcal{F}_{s}^{+}\right.  \right)  =0
\label{rightderivative}%
\end{equation}
$\mathbb{P}_{\mu}$-a.s. for each $t$. Owing to the basic properties of
conditional expectations, this amounts to proving that%
\begin{equation}
\lim_{r\rightarrow t_{+}}\left(  r-t\right)  ^{-1}\mathbb{E}_{\mu}\left(
\mathbb{E}_{\mu}\left(  y_{r}^{\ast}-y_{t}^{\ast}\left\vert \mathcal{F}%
_{t}^{+}\right.  \right)  \left\vert \mathcal{F}_{s}^{+}\right.  \right)  =0
\label{rightderivative2}%
\end{equation}
for each $t$ since $\mathcal{F}_{s}^{+}\subseteq\mathcal{F}_{t}^{+}$.

In order to get (\ref{rightderivative2}) we first show that%
\begin{equation}
\lim_{r\rightarrow t_{+}}\left(  r-t\right)  ^{-1}\mathbb{E}_{\mu}\left(
y_{r}^{\ast}-y_{t}^{\ast}\left\vert \mathcal{F}_{t}^{+}\right.  \right)  =0.
\label{rightderivative3}%
\end{equation}
We have
\begin{align}
&  \mathbb{E}_{\mu}\left(  y_{r}^{\ast}-y_{t}^{\ast}\left\vert \mathcal{F}%
_{t}^{+}\right.  \right) \nonumber\\
&  =\mathbb{E}_{\mu}\left(  \left(  Z_{r}-Z_{t},q\right)  _{\mathbb{R}^{d}%
}\left\vert \mathcal{F}_{t}^{+}\right.  \right)  -\mathbb{E}_{\mu}\left(
\int_{t}^{r}d\tau\left(  b^{\ast}\left(  Z_{\tau},\tau\right)  ,q\right)
_{\mathbb{R}^{d}}\left\vert \mathcal{F}_{t}^{+}\right.  \right)
\label{conditionalexpec}%
\end{align}
according to (\ref{scalarvaluedprocess}), with%
\begin{align}
&  \mathbb{E}_{\mu}\left(  \left(  Z_{r}-Z_{t},q\right)  _{\mathbb{R}^{d}%
}\left\vert \mathcal{F}_{t}^{+}\right.  \right) \label{conditionalexpec1}\\
&  =\mathbb{E}_{\mu}\left(  \left(  Z_{r}-Z_{t},q\right)  _{\mathbb{R}^{d}%
}\left\vert Z_{t}\right.  \right)  =\int_{D}dym^{\ast}\left(  Z_{t}%
,t;y,r\right)  \left(  y-Z_{t},q\right)  _{\mathbb{R}^{d}}%
\end{align}
and%
\begin{align}
&  \mathbb{E}_{\mu}\left(  \int_{t}^{r}d\tau\left(  b^{\ast}\left(  Z_{\tau
},\tau\right)  ,q\right)  _{\mathbb{R}^{d}}\left\vert \mathcal{F}_{t}%
^{+}\right.  \right) \nonumber\\
&  =\int_{t}^{r}d\tau\mathbb{E}_{\mu}\left(  \left(  b^{\ast}\left(  Z_{\tau
},\tau\right)  ,q\right)  _{\mathbb{R}^{d}}\left\vert Z_{t}\right.  \right)
\nonumber\\
&  =\int_{t}^{r}d\tau\int_{D}dym^{\ast}(Z_{t},t;y,\tau)\left(  b^{\ast}\left(
y,\tau\right)  ,q\right)  _{\mathbb{R}^{d}}. \label{conditionalexpec2}%
\end{align}
Now from (\ref{drift2}), (\ref{drift4}) and (\ref{conditionalexpec1})\ we get%
\begin{equation}
\lim_{r\rightarrow t_{+}}\left(  r-t\right)  ^{-1}\mathbb{E}_{\mu}\left(
\left(  Z_{r}-Z_{t},q\right)  _{\mathbb{R}^{d}}\left\vert \mathcal{F}_{t}%
^{+}\right.  \right)  =\left(  b^{\ast}\left(  Z_{t},t\right)  ,q\right)
_{\mathbb{R}^{d}} \label{limit1}%
\end{equation}
$\mathbb{P}_{\mu}$-a.s.. Furthermore, we also claim that%
\begin{equation}
\lim_{r\rightarrow t_{+}}\left(  r-t\right)  ^{-1}\int_{t}^{r}d\tau\int
_{D}dym^{\ast}(Z_{t},t;y,\tau)\left(  b^{\ast}\left(  y,\tau\right)
,q\right)  _{\mathbb{R}^{d}}=\left(  b^{\ast}\left(  Z_{t},t\right)
,q\right)  _{\mathbb{R}^{d}} \label{cruciallimit0}%
\end{equation}
$\mathbb{P}_{\mu}$-a.s. or, equivalently, that%
\begin{equation}
\lim_{r\rightarrow t_{+}}\left(  r-t\right)  ^{-1}\int_{t}^{r}d\tau\int
_{D}dym^{\ast}(Z_{t},t;y,\tau)\left(  b^{\ast}\left(  y,\tau\right)  -b^{\ast
}\left(  Z_{t},t\right)  ,q\right)  _{\mathbb{R}^{d}}=0. \label{cruciallimit}%
\end{equation}
The crucial fact about proving (\ref{cruciallimit}) is that the drift-term
(\ref{drift4}) satisfies the H\"{o}lder continuity estimate%
\[
\left\vert b^{\ast}\left(  y,\tau\right)  -b^{\ast}\left(  Z_{t},t\right)
\right\vert \leq c\left(  \left\vert y-Z_{t}\right\vert ^{\alpha}+\left\vert
\tau-t\right\vert ^{\frac{\alpha}{2}}\right)  ,
\]
which is an easy consequence of the H\"{o}lder properties of the diffusion
matrix $k$ and those of the vector field $l$ stated in Hypotheses (K) and (L),
together with (a)\ of Proposition 1 regarding $v_{\psi}$. Consequently, by
using the very same kind of estimates as we did in the proofs of
(\ref{continuousversion}) and (\ref{vanishingcontribution}) we obtain
successively%
\begin{align*}
0\leq &  \left(  r-t\right)  ^{-1}\int_{t}^{r}d\tau\int_{D}dym^{\ast}%
(Z_{t},t;y,\tau)\left\vert b^{\ast}\left(  y,\tau\right)  -b^{\ast}\left(
Z_{t},t\right)  \right\vert \\
&  \leq c\left(  r-t\right)  ^{-1}\int_{t}^{r}d\tau\left(  \tau-t\right)
^{-\frac{d}{2}}\int_{\mathbb{R}^{d}}dy\exp\left[  -c\frac{\left\vert
y\right\vert ^{2}}{\tau-t}\right]  \left(  \left\vert y\right\vert ^{\alpha
}+\left\vert \tau-t\right\vert ^{\frac{\alpha}{2}}\right) \\
&  \leq c\left(  r-t\right)  ^{-1}\int_{t}^{r}d\tau\left(  \tau-t\right)
^{\frac{\alpha}{2}}\int_{\mathbb{R}^{d}}dy\exp\left[  -c\left\vert
y\right\vert ^{2}\right]  \left\vert y\right\vert ^{\alpha}\\
&  +c\left(  r-t\right)  ^{-1}\int_{t}^{r}d\tau\left(  \tau-t\right)
^{\frac{\alpha}{2}}\int_{\mathbb{R}^{d}}dy\exp\left[  -c\left\vert
y\right\vert ^{2}\right] \\
&  \leq c\left(  r-t\right)  ^{\frac{\alpha}{2}}\rightarrow0
\end{align*}
as $r\rightarrow t_{+}$, which indeed proves (\ref{cruciallimit}). The
combination of (\ref{conditionalexpec}) with (\ref{conditionalexpec2}%
)-(\ref{cruciallimit}) then gives (\ref{rightderivative3}).

It remains to prove that (\ref{rightderivative2}) is a consequence of
(\ref{rightderivative3}). On the one hand we have%
\begin{equation}
\left\vert \int_{D}dym^{\ast}\left(  Z_{t},t;y,r\right)  \left(
y-Z_{t},q\right)  _{\mathbb{R}^{d}}\right\vert \leq c\left(  r-t\right)
\label{uniformbound1}%
\end{equation}
$\mathbb{P}_{\mu}$-a.s. for some positive, non-random and finite constant $c$
since (\ref{limit1}) holds and the right-hand side of (\ref{limit1}) is
uniformly bounded. On the other hand, we also have%
\begin{equation}
\left\vert \int_{t}^{r}d\tau\int_{D}dym^{\ast}(Z_{t},t;y,\tau)\left(  b^{\ast
}\left(  y,\tau\right)  ,q\right)  _{\mathbb{R}^{d}}\right\vert \leq c\left(
r-t\right)  \label{uniformbound2}%
\end{equation}
because of (\ref{cruciallimit0}). Consequently, it follows from
(\ref{conditionalexpec})-(\ref{conditionalexpec2}) and (\ref{uniformbound1}),
(\ref{uniformbound2}) that
\[
\left\vert \mathbb{E}_{\mu}\left(  y_{r}^{\ast}-y_{t}^{\ast}\left\vert
\mathcal{F}_{t}^{+}\right.  \right)  \right\vert \leq c\left(  r-t\right)
\]
$\mathbb{P}_{\mu}$-a.s. for some suitable non-random and finite constant $c$
and all $r,t\in\left[  0,T\right]  $ with $r\geq t$. Therefore,
(\ref{rightderivative3}) indeed implies (\ref{rightderivative2}) by dominated
convergence, so that (\ref{rightderivative}) obtains and thereby the
martingale property (\ref{martingaleproperty}). The proof of the statement
concerning the continuous process (\ref{continuousmartingale2}) follows from
similar arguments, by proving that%
\[
\sup_{t\in\left[  0,T\right]  }\mathbb{E}_{\mu}\left\vert Y_{t}\right\vert
^{2}<+\infty
\]
and that the relation%
\[
\mathbb{E}_{\mu}\left(  y_{s}\left\vert \mathcal{F}_{t}^{-}\right.  \right)
=y_{t}%
\]
holds $\mathbb{P}_{\mu}$-a.s. for all $s\leq$ $t$, where%
\[
y_{t}:=\left(  Z_{t}-Z_{T},q\right)  _{\mathbb{R}^{d}}+\int_{t}^{T}%
d\tau\left(  b\left(  Z_{\tau},\tau\right)  ,q\right)  _{\mathbb{R}^{d}}%
\]
and $\mathcal{F}_{t}^{-}\subseteq\mathcal{F}_{s}^{-}$. \ \ $\blacksquare$

\bigskip

\textsc{Remark.} Let us assume momentarily that $b^{\ast}=b=0$ identically. We
then infer from (\ref{continuousmartingale}) and (\ref{continuousmartingale2})
that $Z_{\tau\in\left[  0,T\right]  }$ is a continuous, square-integrable
martingale for both filtrations $\mathcal{F}_{\tau\in\left[  0,T\right]  }%
^{+}$ and $\mathcal{F}_{\tau\in\left[  0,T\right]  }^{-}$ simultaneously.
Consequently we have%
\[
\mathbb{E}_{\mu}\left(  \left(  Z_{t}-Z_{s},q\right)  _{\mathbb{R}^{d}}%
^{2}\left\vert \mathcal{F}_{s}^{+}\right.  \right)  =\mathbb{E}_{\mu}\left(
\left(  Z_{t},q\right)  _{\mathbb{R}^{d}}^{2}\left\vert \mathcal{F}_{s}%
^{+}\right.  \right)  -{}\left(  Z_{s},q\right)  _{\mathbb{R}^{d}}^{2}%
\]
$\mathbb{P}_{\mu}$-a.s. for all $t\geq s$ and every $q\in\mathbb{R}^{d}$, and
at the same time%
\[
\mathbb{E}_{\mu}\left(  \left(  Z_{t}-Z_{s},q\right)  _{\mathbb{R}^{d}}%
^{2}\left\vert \mathcal{F}_{t}^{-}\right.  \right)  =\mathbb{E}_{\mu}\left(
\left(  Z_{s},q\right)  _{\mathbb{R}^{d}}^{2}\left\vert \mathcal{F}_{t}%
^{-}\right.  \right)  -{}\left(  Z_{t},q\right)  _{\mathbb{R}^{d}}^{2}.
\]
By averaging both expressions we then obtain%
\begin{align*}
&  \mathbb{E}_{\mu}\left(  Z_{t}-Z_{s},q\right)  _{\mathbb{R}^{d}}^{2}\\
&  =\mathbb{E}_{\mu}\left(  Z_{t},q\right)  _{\mathbb{R}^{d}}^{2}%
-\mathbb{E}_{\mu}\left(  Z_{s},q\right)  _{\mathbb{R}^{d}}^{2}\\
&  =\mathbb{E}_{\mu}\left(  Z_{s},q\right)  _{\mathbb{R}^{d}}^{2}%
-\mathbb{E}_{\mu}\left(  Z_{t},q\right)  _{\mathbb{R}^{d}}^{2}=0
\end{align*}
for every $q\in\mathbb{R}^{d}$, so that $Z_{t}=Z_{s}$ $\mathbb{P}_{\mu}$-a.s.
for all $t\geq s$. Conversely, if the function $\tau\mapsto Z_{\tau}$ is
constant $\mathbb{P}_{\mu}$-a.s. we obviously have $b^{\ast}=b=0$ identically.
From now on we will exclude this exceptional case from our considerations for
reasons that will become apparent from the statement of Theorem 3 below and
its proof.

\bigskip

It is worth emphasizing the fact that the existence of the Wiener processes we
are looking for is conditioned by that of the reversible Markov diffusion
$Z_{\tau\in\left[  0,T\right]  }$ and not the other way around, as there is
absolutely nothing stochastic in Problems (\ref{parabolicproblem}) and
(\ref{adjointproblem}). Furthermore, from a technical point of view all the
stochastic integrals defined below involve continuous square-integrable
martingales as integrators, as for instance in Chapter 3 of
\cite{karatzasshreve}. Then, the precise result about $Z_{\tau\in\left[
0,T\right]  }$ being a reversible It\^{o} diffusion is the following, in which
we write $k^{\frac{1}{2}}\left(  x,t\right)  $ for the positive square root of
the diffusion matrix $k(x,t)$ and assume that $\tau\mapsto Z_{\tau}$ is not
constant $\mathbb{P}_{\mu}$-a.s. according to the preceding remark.

\bigskip

\textbf{Theorem 3. }\textit{Assume that the hypotheses of Lemmas 1 and 2 are
valid. Assume furthermore that }$\mu$ \textit{is given by
(\ref{probabilitymeasure}) and (\ref{specialclass}), and let }$Z_{\tau
\in\left[  0,T\right]  }$ \textit{be the reversible Markov diffusion of
Theorem 2. Then the following statements hold:}

\textit{(a) There exists a }$d$\textit{-dimensional Wiener process }%
$W_{\tau\in\left[  0,T\right]  }^{\ast}$ \textit{such that the relation}%
\begin{equation}
Z_{t}=Z_{0}+\int_{0}^{t}d\tau b^{\ast}\left(  Z_{\tau},\tau\right)  +\int
_{0}^{t}k^{\frac{1}{2}}\left(  Z_{\tau},\tau\right)  d^{+}W_{\tau}^{\ast}
\label{itoprocess1}%
\end{equation}
\textit{holds }$\mathbb{P}_{\mu}$\textit{-a.s. for every} $t\in\left[
0,T\right]  $\textit{, with }$b^{\ast}$ \textit{given by (\ref{drift4}) and
the forward stochastic integral defined with respect to }$\mathcal{F}_{\tau
\in\left[  0,T\right]  }^{+}$.

\textit{(b) There exists a }$d$-\textit{dimensional Wiener process }%
$W_{\tau\in\left[  0,T\right]  }$ \textit{such that the relation}%
\begin{equation}
Z_{t}=Z_{T}-\int_{t}^{T}d\tau b\left(  Z_{\tau},\tau\right)  -\int_{t}%
^{T}k^{\frac{1}{2}}\left(  Z_{\tau},\tau\right)  d^{-}W_{\tau}
\label{itoprocess2}%
\end{equation}
\textit{holds }$\mathbb{P}_{\mu}$\textit{-a.s. for every} $t\in\left[
0,T\right]  $\textit{, with }$b$\textit{ given by (\ref{drift5}) and the
backward stochastic integral defined with respect to }$\mathcal{F}_{\tau
\in\left[  0,T\right]  }^{-}$.

\bigskip

\textbf{Proof. }We define $W_{\tau\in\left[  0,T\right]  }^{\ast}$ by%
\begin{equation}
W_{t}^{\ast}:=\int_{0}^{t}k^{-\frac{1}{2}}\left(  Z_{\tau},\tau\right)
d^{+}Y_{\tau}^{\ast} \label{wienerprocess1}%
\end{equation}
with respect to the continuous square-integrable martingale $Y_{\tau\in\left[
0,T\right]  }^{\ast}$ given by (\ref{continuousmartingale}) and $\mathcal{F}%
_{\tau\in\left[  0,T\right]  }^{+}$. In order to show that
(\ref{wienerprocess1}) makes sense and indeed defines a Wiener process, we
begin by proving that the quadratic variation of the martingale $y_{\tau
\in\left[  0,T\right]  }^{\ast}$ is the absolutely continuous process%
\begin{equation}
\left\langle y^{\ast}\right\rangle _{t}:=\int_{0}^{t}d\tau\left(  k\left(
Z_{\tau},\tau\right)  q,q\right)  _{\mathbb{R}^{d}} \label{quadraticvariation}%
\end{equation}
where $q\in\mathbb{R}^{d}$. According to the Doob-Meyer decomposition, this
amounts to showing that the process defined by%
\begin{equation}
z_{t}^{\ast}:=\left(  y_{t}^{\ast}\right)  ^{2}-\int_{0}^{t}d\tau\left(
k\left(  Z_{\tau},\tau\right)  q,q\right)  _{\mathbb{R}^{d}} \label{doobmeyer}%
\end{equation}
is a martingale relative to $\mathcal{F}_{\tau\in\left[  0,T\right]  }^{+}$.
In order to achieve this we proceed as we did in the proof of Lemma 4, by
proving that the relation%
\begin{align}
&  \lim_{r\rightarrow t_{+}}\left(  r-t\right)  ^{-1}\mathbb{E}_{\mu}\left(
z_{r}^{\ast}-z_{t}^{\ast}\left\vert \mathcal{F}_{s}^{+}\right.  \right)
\nonumber\\
&  =\lim_{r\rightarrow t_{+}}\left(  r-t\right)  ^{-1}\mathbb{E}_{\mu}\left(
\mathbb{E}_{\mu}\left(  z_{r}^{\ast}-z_{t}^{\ast}\left\vert \mathcal{F}%
_{t}^{+}\right.  \right)  \left\vert \mathcal{F}_{s}^{+}\right.  \right)  =0
\label{righthanderivative}%
\end{align}
holds, where for $0\leq t<r\leq T$ the inner conditional expectation is given
by
\begin{align}
&  \mathbb{E}_{\mu}\left(  z_{r}^{\ast}-z_{t}^{\ast}\left\vert \mathcal{F}%
_{t}^{+}\right.  \right) \nonumber\\
&  =\mathbb{E}_{\mu}\left(  \left(  y_{r}^{\ast}-y_{t}^{\ast}\right)
^{2}\left\vert \mathcal{F}_{t}^{+}\right.  \right)  -\int_{t}^{r}%
d\tau\mathbb{E}_{\mu}\left(  \left(  k\left(  Z_{\tau},\tau\right)
q,q\right)  _{\mathbb{R}^{d}}\left\vert \mathcal{F}_{t}^{+}\right.  \right)
\label{innerconditionalexpec}%
\end{align}
since $y_{\tau\in\left[  0,T\right]  }^{\ast}$ is a martingale relative to the
filtration $\mathcal{F}_{\tau\in\left[  0,T\right]  }^{+}$.

We first show that%

\begin{equation}
\lim_{r\rightarrow t_{+}}\left(  r-t\right)  ^{-1}\mathbb{E}_{\mu}\left(
\left(  y_{r}^{\ast}-y_{t}^{\ast}\right)  ^{2}\left\vert \mathcal{F}_{t}%
^{+}\right.  \right)  =\left(  k\left(  Z_{t},t\right)  q,q\right)
_{\mathbb{R}^{d}} \label{conditionalexpec3}%
\end{equation}
$\mathbb{P}_{\mu}$-a.s. for every $t$. On the one hand, from
(\ref{scalarvaluedprocess}) we have%
\begin{align}
&  \left(  y_{r}^{\ast}-y_{t}^{\ast}\right)  ^{2}\nonumber\\
&  =\left(  Z_{r}-Z_{t},q\right)  _{\mathbb{R}^{d}}^{2}+\left(  \int_{t}%
^{r}d\tau\left(  b^{\ast}(Z_{\tau},\tau\right)  ,q)_{\mathbb{R}^{d}}\right)
^{2}\nonumber\\
&  -2\left(  Z_{r}-Z_{t},q\right)  _{\mathbb{R}^{d}}\int_{t}^{r}d\tau\left(
b^{\ast}(Z_{\tau},\tau\right)  ,q)_{\mathbb{R}^{d}} \label{squaredprocess}%
\end{align}
and we note that%
\begin{align*}
&  \lim_{r\rightarrow t_{+}}\left(  r-t\right)  ^{-1}\mathbb{E}_{\mu}\left(
\left(  Z_{r}-Z_{t},q\right)  _{\mathbb{R}^{d}}^{2}\left\vert \mathcal{F}%
_{t}^{+}\right.  \right) \\
&  =\lim_{r\rightarrow t_{+}}\left(  r-t\right)  ^{-1}\mathbb{E}_{\mu}\left(
\left(  Z_{r}-Z_{t},q\right)  _{\mathbb{R}^{d}}^{2}\left\vert Z_{t}\right.
\right)  =\left(  k\left(  Z_{t},t\right)  q,q\right)  _{\mathbb{R}^{d}}%
\end{align*}
since $\tau\mapsto Z_{\tau}$ is not constant $\mathbb{P}_{\mu}$-a.s.. Indeed
this follows immediately from (\ref{diffusion}), which implies the relation%
\[
\lim_{r\rightarrow t_{+}}\left(  r-t\right)  ^{-1}\int_{D}dym^{\ast}%
(Z_{t},t;y,r)\left(  y-Z_{t},q\right)  _{\mathbb{R}^{d}}^{2}=\left(  k\left(
Z_{t},t\right)  q,q\right)  _{\mathbb{R}^{d}}%
\]
by switching to the quadratic form formulation. On the other hand, by
dominated convergence the remaining terms on the right-hand side of
(\ref{squaredprocess}) do not contribute to the conditional expectation
(\ref{conditionalexpec3}), for%
\[
\lim_{r\rightarrow t_{+}}\left(  r-t\right)  ^{-1}\left(  \int_{t}^{r}%
d\tau\left(  b^{\ast}(Z_{\tau},\tau\right)  ,q)_{\mathbb{R}^{d}}\right)
^{2}=0
\]
and%
\begin{align*}
&  \lim_{r\rightarrow t_{+}}\left(  r-t\right)  ^{-1}\left(  Z_{r}%
-Z_{t},q\right)  _{\mathbb{R}^{d}}\int_{t}^{r}d\tau\left(  b^{\ast}(Z_{\tau
},\tau\right)  ,q)_{\mathbb{R}^{d}}\\
&  =\left(  b^{\ast}(Z_{t},t\right)  ,q)_{\mathbb{R}^{d}}\lim_{r\rightarrow
t_{+}}\left(  Z_{r}-Z_{t},q\right)  _{\mathbb{R}^{d}}=0
\end{align*}
$\mathbb{P}_{\mu}$-a.s. for every $t$, as a consequence of the boundedness and
the continuity of $b^{\ast}$ and $Z_{\tau\in\left[  0,T\right]  }$. Therefore
(\ref{conditionalexpec3}) holds.

Next, we observe that we also have%
\begin{equation}
\lim_{r\rightarrow t_{+}}\left(  r-t\right)  ^{-1}\int_{t}^{r}d\tau
\mathbb{E}_{\mu}\left(  \left(  k\left(  Z_{\tau},\tau\right)  q,q\right)
_{\mathbb{R}^{d}}\left\vert \mathcal{F}_{t}^{+}\right.  \right)  =\left(
k\left(  Z_{t},t\right)  q,q\right)  _{\mathbb{R}^{d}}
\label{conditionalexpec4}%
\end{equation}
$\mathbb{P}_{\mu}$-a.s. for every $t$ or, equivalently, that%
\[
\lim_{r\rightarrow t_{+}}\left(  r-t\right)  ^{-1}\int_{t}^{r}d\tau\int
_{D}dym^{\ast}(Z_{t},t;y,\tau)\left(  (k\left(  y,\tau\right)  -k\left(
Z_{t},t\right)  )q,q\right)  _{\mathbb{R}^{d}}=0.
\]
Indeed this relation follows from exactly the same arguments as those which
led to (\ref{cruciallimit}), since the matrix elements of $k$ are all jointly
H\"{o}lder continuous relative to space-time variables according to Hypothesis
(K). Relations (\ref{innerconditionalexpec}), (\ref{conditionalexpec3}) and
(\ref{conditionalexpec4}) then imply%
\[
\lim_{r\rightarrow t_{+}}\left(  r-t\right)  ^{-1}\mathbb{E}_{\mu}\left(
z_{r}^{\ast}-z_{t}^{\ast}\left\vert \mathcal{F}_{t}^{+}\right.  \right)  =0
\]
$\mathbb{P}_{\mu}$-a.s. for every $t$, so that (\ref{righthanderivative})
follows from a dominated convergence argument similar to that given in the
proof of Lemma 4. Thus $z_{\tau\in\left[  0,T\right]  }^{\ast}$ is a
martingale relative to $\mathcal{F}_{\tau\in\left[  0,T\right]  }^{+}$, and
consequently (\ref{quadraticvariation}) is indeed the quadratic variation
process of (\ref{scalarvaluedprocess}).

This implies (\ref{wienerprocess1}) makes sense in that $k^{-\frac{1}{2}%
}\left(  Z_{\tau},\tau\right)  $ is an admissible integrand there, and defines
a continuous, square-integrable martingale whose quadratic variation is given
by%
\[
\left\langle \left(  W^{\ast},q\right)  _{\mathbb{R}^{d}}\right\rangle
_{t}=t\left\vert q\right\vert ^{2}%
\]
for each $t\in\left[  0,T\right]  $ and every $q\in\mathbb{R}^{d}$. This
proves that $W_{\tau\in\left[  0,T\right]  }^{\ast}$ is a Wiener process, and
furthermore that the combination of (\ref{wienerprocess1}) with
(\ref{continuousmartingale}) leads to%
\begin{align*}
&  \int_{0}^{t}k^{\frac{1}{2}}\left(  Z_{\tau},\tau\right)  d^{+}W_{\tau
}^{\ast}\\
&  =\int_{0}^{t}d^{+}Y_{\tau}^{\ast}=Z_{t}-Z_{0}-\int_{0}^{t}d\tau b^{\ast
}\left(  Z_{\tau},\tau\right)
\end{align*}
\textit{ }$\mathbb{P}_{\mu}$-a.s. for every $t\in\left[  0,T\right]  $, which
is (\ref{itoprocess1}). The proof of (\ref{itoprocess2}) with
\begin{equation}
W_{t}:=-\int_{t}^{T}k^{-\frac{1}{2}}\left(  Z_{\tau},\tau\right)  d^{-}%
Y_{\tau} \label{wienerprocess2}%
\end{equation}
defined with respect to (\ref{continuousmartingale2}) is similar, with this
time%
\[
z_{t}:=\left(  y_{t}\right)  ^{2}-\int_{t}^{T}d\tau\left(  k\left(  Z_{\tau
},\tau\right)  q,q\right)  _{\mathbb{R}^{d}}%
\]
a martingale with respect to the decreasing filtration $\mathcal{F}_{\tau
\in\left[  0,T\right]  }^{-}$. \ \ $\blacksquare$

\bigskip

We have already noted that (\ref{bernsteindensity}) provides important
information regarding $Z_{\tau\in\left[  0,T\right]  }$ in terms of the
solutions to (\ref{parabolicproblem}) and (\ref{adjointproblem}). It is,
therefore, natural to ask whether those solutions can, in turn, be represented
as suitable expectations of some functionals of $Z_{\tau\in\left[  0,T\right]
}$. This is indeed possible but not \textit{a priori} evident since the
elliptic operators on the right-hand side of (\ref{parabolicproblem}) and
(\ref{adjointproblem}) are\textit{ not} the generators of $Z_{\tau\in\left[
0,T\right]  }$. The problem lies, of course, in the presence of the
logarithmic terms in (\ref{drift4}), (\ref{drift5}), and in the next result we
show how to do do away with them by means of suitable Girsanov transformations.

\bigskip

\textbf{Corollary.} \textit{Assume that the same hypotheses as in Theorem 3
are valid, and let }$Z_{\tau\in\left[  0,T\right]  }$ \textit{be the
reversible Markov diffusion of Theorem 2. Then, aside from }$\mathbb{P}_{\mu}$
\textit{there exist two probability measures }$\mathbb{P}_{\mu}^{\pm}$\textit{
on} $\left(  \Omega,\mathcal{F}\right)  $ \textit{such that} \textit{the
following statements hold:}

\textit{(a) There exists a }$d$\textit{-dimensional Wiener process
}$\widetilde{W}_{\tau\in\left[  0,T\right]  }^{\ast}$ \textit{relative to
}$\mathcal{F}_{\tau\in\left[  0,T\right]  }^{+}$ \textit{such that the
relation }%
\begin{equation}
Z_{t}=Z_{0}+\int_{0}^{t}d\tau a^{\ast}\left(  Z_{\tau},\tau\right)  +\int
_{0}^{t}k^{\frac{1}{2}}\left(  Z_{\tau},\tau\right)  d^{+}\widetilde{W}_{\tau
}^{\ast} \label{itoprocess3}%
\end{equation}
\textit{holds }$\mathbb{P}_{\mu}^{+}$\textit{-a.s. for every} $t\in\left[
0,T\right]  $\textit{, with }$a^{\ast}$ \textit{given by (\ref{driftterm1}).}

\textit{(b) There exists a }$d$-\textit{dimensional Wiener process
}$\widetilde{W}_{\tau\in\left[  0,T\right]  }$ relative to $\mathcal{F}%
_{\tau\in\left[  0,T\right]  }^{-}$ \textit{such that the relation}%
\begin{equation}
Z_{t}=Z_{T}-\int_{t}^{T}d\tau a\left(  Z_{\tau},\tau\right)  -\int_{t}%
^{T}k^{\frac{1}{2}}\left(  Z_{\tau},\tau\right)  d^{-}\widetilde{W}_{\tau}
\label{itoprocess4}%
\end{equation}
\textit{holds }$\mathbb{P}_{\mu}^{-}$\textit{-a.s. for every} $t\in\left[
0,T\right]  $\textit{, with }$a$ \textit{given by (\ref{driftterm2}).}

\bigskip

\textbf{Proof.} Let us define%
\begin{equation}
\widetilde{W}_{t}^{\ast}:=\int_{0}^{t}d\tau X_{\tau}+W_{t}^{\ast}
\label{wienerprocess3}%
\end{equation}
for every $t\in\left[  0,T\right]  $, where%
\begin{equation}
X_{t}:=k^{\frac{1}{2}}\left(  Z_{t},t\right)  \nabla_{x}\ln v_{\psi}(Z_{t},t).
\label{auxiliaryprocess}%
\end{equation}
From the hypotheses regarding $k$ and the properties of $v_{\psi}$ and
$Z_{\tau\in\left[  0,T\right]  }$, it is clear that (\ref{auxiliaryprocess})
defines a continuous process adapted to the filtration $\mathcal{F}_{\tau
\in\left[  0,T\right]  }^{+}$. It is also bounded $\mathbb{P}_{\mu}$-a.s. by a
non-random positive constant, so that we have%
\[
\mathbb{E}_{\mu}\exp\left[  \frac{1}{2}\int_{0}^{T}d\tau\left\vert X_{\tau
}\right\vert ^{2}\right]  <+\infty
\]
and thereby%
\[
\mathbb{E}_{\mu}\exp\left[  -\int_{0}^{T}\left(  X_{\tau},d^{+}W_{\tau}^{\ast
}\right)  _{\mathbb{R}^{d}}-\frac{1}{2}\int_{0}^{T}d\tau\left\vert X_{\tau
}\right\vert ^{2}\right]  =1
\]
according to Theorem 12 in Section 1 of Chapter 3 of \cite{gihmanskohorod}.
For every $E\in\mathcal{F}$ we then define%
\[
\mathbb{P}_{\mu}^{+}(E):=\int_{E}d\mathbb{P}_{\mu}\exp\left[  -\int_{0}%
^{T}\left(  X_{\tau},d^{+}W_{\tau}^{\ast}\right)  _{\mathbb{R}^{d}}-\frac
{1}{2}\int_{0}^{T}d\tau\left\vert X_{\tau}\right\vert ^{2}\right]  \mathit{,}%
\]
which indeed makes (\ref{wienerprocess3}) a Wiener process on $(\Omega
,\mathcal{F},\mathbb{P}_{\mu}^{+})$ relative to $\mathcal{F}_{\tau\in\left[
0,T\right]  }^{+}$ according to Girsanov's standard construction (see, for
instance, \cite{gihmanskohorod}). Furthermore we have%
\begin{align*}
&  \int_{0}^{t}k^{\frac{1}{2}}\left(  Z_{\tau},\tau\right)  d^{+}\widetilde
{W}_{\tau}^{\ast}\\
&  =\int_{0}^{t}d\tau k\left(  Z_{\tau},\tau\right)  \nabla_{x}\ln v_{\psi
}(Z_{\tau},\tau)+\int_{0}^{t}k^{\frac{1}{2}}\left(  Z_{\tau},\tau\right)
d^{+}W_{\tau}^{\ast}%
\end{align*}
$\mathbb{P}_{\mu}$-a.s. and $\mathbb{P}_{\mu}^{+}$-a.s. for every $t\in\left[
0,T\right]  $, so that (\ref{itoprocess1}) reduces to (\ref{itoprocess3}). The
proof of (\ref{itoprocess4}) with%
\[
\widetilde{W}_{t}:=\int_{t}^{T}d\tau k^{\frac{1}{2}}\left(  Z_{\tau}%
,\tau\right)  \nabla_{x}\ln u_{\varphi}(Z_{\tau},\tau)+W_{t}%
\]
is similar and therefore omitted. \ \ $\blacksquare$

\bigskip

Let us now write $\mathbb{E}_{\mu}^{\pm}$ for the expectation functional on
$(\Omega,\mathcal{F},\mathbb{P}_{\mu}^{\pm})$, and by $\mathbb{E}_{\mu
,x,t}^{\pm}$ the conditional expectations corresponding to setting $Z_{t}=x$
for an arbitrary $(x,t)\in D\times\left[  0,T\right]  $. We then get the
following Feynman-Kac representations for the solutions to
(\ref{parabolicproblem}) and (\ref{adjointproblem}) we alluded to above.

\bigskip

\textbf{Theorem 4.} \textit{Assume that the same hypotheses as in the
Corollary are valid, and let }$Z_{\tau\in\left[  0,T\right]  }$ \textit{be the
reversible Markov diffusion of that Corollary. Then the following statements
hold:}

\textit{(a) The unique classical positive solution to (\ref{parabolicproblem})
may be written as}%
\begin{equation}
u_{\varphi}(x,t)=\mathbb{E}_{\mu,x,t}^{-}\left(  \exp\left[  -\int_{0}%
^{t}d\sigma V\left(  Z_{\sigma},\sigma\right)  \right]  \varphi\left(
Z_{0}\right)  \right)  \label{feynmankac2}%
\end{equation}
\textit{for all} $(x,t)\in D\times\left[  0,T\right]  .$

\textit{(b) If in addition the vector-field }$l$ \textit{satisfies}%
\begin{equation}
\operatorname{div}_{x}l(x,t)=0 \label{divergencecondition}%
\end{equation}
\textit{for each }$t\in\left[  0,T\right]  $\textit{, then the unique
classical positive solution to (\ref{adjointproblem}) may be written as}%
\begin{equation}
v_{\psi}(x,t)=\mathbb{E}_{\mu,x,t}^{+}\left(  \exp\left[  -\int_{t}^{T}d\sigma
V\left(  Z_{\sigma},\sigma\right)  \right]  \psi\left(  Z_{T}\right)  \right)
\label{feynmankac1}%
\end{equation}
\textit{for all} $(x,t)\in D\times\left[  0,T\right]  .$

\bigskip

\textbf{Proof.} Regarding the proof of (a) we have%
\begin{align}
&  \frac{1}{2}\operatorname{div}_{x}\left(  k(x,t)\nabla_{x}u_{\varphi
}(x,t)\right)  -\left(  l(x,t),\nabla_{x}u_{\varphi}(x,t)\right)
_{\mathbb{R}^{d}}\nonumber\\
&  =\frac{1}{2}\left(  k(x,t)\nabla_{x},\nabla_{x}u_{\varphi}(x,t)\right)
_{\mathbb{R}^{d}}-\left(  a(x,t),\nabla_{x}u_{\varphi}(x,t)\right)
_{\mathbb{R}^{d}} \label{backwardgenerator}%
\end{align}
by using (\ref{driftterm2}), so that the sum of the first two terms on the
right-hand side of the first line of (\ref{parabolicproblem}) identifies with
the generator of the backward diffusion (\ref{itoprocess4}). We then apply the
corresponding \textit{backward} It\^{o} formula%
\begin{align}
&  F\left(  Z_{t},t\right) \nonumber\\
&  =F\left(  Z_{s},s\right)  +\int_{s}^{t}d\tau\frac{\partial F}{\partial
s}\left(  Z_{\tau},\tau\right)  -\frac{1}{2}\int_{s}^{t}d\tau\left(
k(Z_{\tau},\tau)\nabla_{x},\nabla_{x}F(Z_{\tau},\tau)\right)  _{\mathbb{R}%
^{d}}\nonumber\\
&  +\int_{s}^{t}d\tau\left(  a(Z_{\tau},\tau),\nabla_{x}F(Z_{\tau}%
,\tau)\right)  _{\mathbb{R}^{d}}+\int_{s}^{t}\left(  \nabla_{x}F(Z_{\tau}%
,\tau),k^{\frac{1}{2}}\left(  Z_{\tau},\tau\right)  d^{-}\widetilde{W}_{\tau
}\right)  _{\mathbb{R}^{d}} \label{backwarditoformula}%
\end{align}
to the function%
\[
F\left(  x,s\right)  :=\exp\left[  -\int_{s}^{t}d\sigma V\left(  Z_{\sigma
},\sigma\right)  \right]  u_{\varphi}(x,s)
\]
for each $x\in\overline{D}$ and all $s,t\in\left[  0,T\right]  $ satisfying
$t\geq s$, by noticing in particular the crucial minus sign of the last term
of the second line to (\ref{backwarditoformula}). By substituting%
\[
\frac{\partial F}{\partial s}\left(  x,s\right)  =\exp\left[  -\int_{s}%
^{t}d\sigma V\left(  Z_{\sigma},\sigma\right)  \right]  \left(  \partial
_{s}u_{\varphi}\left(  x,s\right)  +V(Z_{s},s)u_{\varphi}(x,s)\right)
\]
and%
\[
\nabla_{x}F(x,s)=\exp\left[  -\int_{s}^{t}d\sigma V\left(  Z_{\sigma}%
,\sigma\right)  \right]  \nabla_{x}u_{\varphi}(x,s)
\]
into (\ref{backwarditoformula}), and by using (\ref{backwardgenerator})
together with the first relation in (\ref{parabolicproblem}), we obtain%
\begin{align}
&  u_{\varphi}\left(  Z_{t},t\right) \nonumber\\
&  =\exp\left[  -\int_{s}^{t}d\sigma V\left(  Z_{\sigma},\sigma\right)
\right]  u_{\varphi}\left(  Z_{s},s\right) \nonumber\\
&  +\int_{s}^{t}\exp\left[  -\int_{\tau}^{t}d\sigma V\left(  Z_{\sigma}%
,\sigma\right)  \right]  \left(  \nabla_{x}u_{\varphi}(Z_{\tau},\tau
),k^{\frac{1}{2}}\left(  Z_{\tau},\tau\right)  d^{-}\widetilde{W}_{\tau
}\right)  _{\mathbb{R}^{d}} \label{stochasticrepresentation}%
\end{align}
where the stochastic integral in (\ref{stochasticrepresentation}) is finite
$\mathbb{P}_{\mu}^{-}$-a.s.. The desired result (\ref{feynmankac2}) then
follows by setting $s=0$ in the preceding relation, by using the initial
condition in (\ref{parabolicproblem}) and by taking the conditional
expectation $\mathbb{E}_{\mu,x,t}^{-}$ of the resulting equality.

As for the proof of (b), we first observe that (\ref{adjointproblem}) may be
rewritten as%
\begin{align}
-\partial_{t}v(x,t)  &  =\frac{1}{2}\operatorname{div}_{x}\left(
k(x,t)\nabla_{x}v(x,t)\right)  +\left(  l(x,t),\nabla_{x}v(x,t)\right)
_{\mathbb{R}^{d}}-V(x,t)v(x,t),\nonumber\\
(x,t)  &  \in D\times\left[  0,T\right)  ,\nonumber\\
v(x,T)  &  =\psi(x),\text{ \ \ }x\in D,\nonumber\\
\frac{\partial v(x,t)}{\partial n_{k}(x,t)}  &  =0,\text{ \ \ \ }%
(x,t)\in\partial D\times\left[  0,T\right)  , \label{adjointproblembis}%
\end{align}
as a consequence of (\ref{divergencecondition}). Furthermore we have%
\begin{align*}
&  \frac{1}{2}\operatorname{div}_{x}\left(  k(x,t)\nabla_{x}v_{\varphi
}(x,t)\right)  +\left(  l(x,t),\nabla_{x}v_{\varphi}(x,t)\right)
_{\mathbb{R}^{d}}\\
&  =\frac{1}{2}\left(  k(x,t)\nabla_{x},\nabla_{x}v_{\varphi}(x,t)\right)
_{\mathbb{R}^{d}}+\left(  a^{\ast}(x,t),\nabla_{x}v_{\varphi}(x,t)\right)
_{\mathbb{R}^{d}}%
\end{align*}
by using (\ref{driftterm1}), so that the sum of the first two terms on the
right-hand side of the first line of (\ref{adjointproblembis}) identifies this
time with the generator of the forward diffusion (\ref{itoprocess3}). The
remaining part of the argument then consists in applying the usual
\textit{forward} It\^{o} formula%
\begin{align}
&  F\left(  Z_{t},t\right) \nonumber\\
&  =F\left(  Z_{s},s\right)  +\int_{s}^{t}d\tau\frac{\partial F}{\partial
t}\left(  Z_{\tau},\tau\right)  +\frac{1}{2}\int_{s}^{t}d\tau\left(
k(Z_{\tau},\tau)\nabla_{x},\nabla_{x}F(Z_{\tau},\tau)\right)  _{\mathbb{R}%
^{d}}\nonumber\\
&  +\int_{s}^{t}d\tau\left(  a^{\ast}(Z_{\tau},\tau),\nabla_{x}F(Z_{\tau}%
,\tau)\right)  _{\mathbb{R}^{d}}+\int_{s}^{t}\left(  \nabla_{x}F(Z_{\tau}%
,\tau),k^{\frac{1}{2}}\left(  Z_{\tau},\tau\right)  d^{+}\widetilde{W}%
_{t}^{\ast}\right)  _{\mathbb{R}^{d}} \label{forwarditoformula}%
\end{align}
to the function given by%

\[
F\left(  x,t\right)  :=\exp\left[  -\int_{s}^{t}d\sigma V\left(  Z_{\sigma
},\sigma\right)  \right]  v_{\psi}(x,t)
\]
for every $\left(  x,t\right)  \in\overline{D}\times\left[  s,T\right]  $,
which eventually allows us to proceed as in the first part of the proof by
using (\ref{adjointproblembis}) to obtain (\ref{feynmankac1}).
\ \ $\blacksquare$

\bigskip

\textsc{Remarks.} (1) The backward It\^{o} formula (\ref{backwarditoformula})
with the appropriate sign change has its origins in the theory developed in
Section 13 of \cite{nelson1}, particularly in Relation (2) of that section
which in fact refers to a simpler situation. Its proof is similar to that of
its forward counterpart (\ref{forwarditoformula}) which is stated in many
places in a more general form, for instance in Theorem 3.6 of Chapter 3 in
\cite{karatzasshreve} or in Theorem 5.1 of Chapter 2 in \cite{ikedawatanabe}.
In this connection it is worth mentioning that we are not aware of any other
\textit{direct} proofs of Feynman-Kac representations such as
(\ref{feynmankac2}) regarding the solutions to \textit{forward }non-autonomous
equations of the form (\ref{parabolicproblem}). Such a proof is indeed made
possible in our case thanks to the \textit{backward} It\^{o} equation
(\ref{itoprocess4}), which allows us to carry out the proofs of
(\ref{feynmankac2}) and (\ref{feynmankac1}) quite independently as there are
two It\^{o} equations for the same process at our disposal.

(2) It is worth focusing again on two important aspects of our constructions,
namely, on the one hand (\ref{bernsteindensity}) which provides information on
$Z_{\tau\in\left[  0,T\right]  }$ from the knowledge of (\ref{forwardsolution}%
) and (\ref{backwardsolution}) and, on the other hand, the representations
(\ref{feynmankac2}) and (\ref{feynmankac1}) that make those two solutions
emerge as functionals of $Z_{\tau\in\left[  0,T\right]  }$. In the latter case
we remark that the comparison of (\ref{feynmankac2}) and (\ref{feynmankac1})
with (\ref{forwardsolution}) and (\ref{backwardsolution}) gives%
\begin{align*}
\mathbb{E}_{\mu,x,t}^{-}\left(  \exp\left[  -\int_{0}^{t}d\sigma V\left(
Z_{\sigma},\sigma\right)  \right]  \varphi\left(  Z_{0}\right)  \right)   &
=\int_{D}dyg_{A}(x,t;y,0)\varphi(y),\\
\mathbb{E}_{\mu,x,t}^{+}\left(  \exp\left[  -\int_{t}^{T}d\sigma V\left(
Z_{\sigma},\sigma\right)  \right]  \psi\left(  Z_{T}\right)  \right)   &
=\int_{D}dyg_{A}^{\ast}(x,t;y,T)\psi(y),
\end{align*}
respectively, that is, a concrete representation of the expectation
functionals $\mathbb{E}_{\mu,x,t}^{\pm}$. In particular, for $V=0$ we obtain%
\begin{align*}
\mathbb{E}_{\mu,x,t}^{-}\varphi\left(  Z_{0}\right)   &  =\int_{D}%
dyg_{V=0}(x,t;y,0)\varphi(y),\\
\mathbb{E}_{\mu,x,t}^{+}\psi\left(  Z_{T}\right)   &  =\int_{D}dyg_{V=0}%
^{\ast}(x,t;y,T)\psi(y),
\end{align*}
respectively, with $g_{V=0}$ and $g_{V=0}^{\ast}$ the Green functions
associated with (\ref{parabolicproblem}) and (\ref{adjointproblem}) in this case.

\bigskip

In the next section we illustrate some of the preceding results.

\section{Two examples}

In the first example and for the sake of clarity we investigate the case of
the simplest possible heat equation and its adjoint on a one-dimensional domain.

\bigskip

\textsc{Example 1}. Let us consider the initial-boundary value problem%
\begin{align}
\partial_{t}u(x,t)  &  =\frac{1}{2}\partial_{xx}u(x,t),\text{ \ \ }%
(x,t)\in\left(  0,1\right)  \times\left(  0,T\right]  ,\text{\ }\nonumber\\
u(x,0)  &  =\varphi(x),\text{ \ \ }x\in\left(  0,1\right)  ,\nonumber\\
\partial_{x}u(0,t)  &  =\partial_{x}u(1,t)=0,\text{ \ \ \ }t\in\left(
0,T\right]  \label{example1forward}%
\end{align}
and the corresponding final-boundary value problem for the adjoint equation%
\begin{align}
-\partial_{t}v(x,t)  &  =\frac{1}{2}\partial_{xx}v(x,t),\text{ \ \ }%
(x,t)\in\left(  0,1\right)  \times\left[  0,T\right)  ,\nonumber\\
v(x,T)  &  =\psi(x),\text{ \ \ }x\in\left(  0,1\right)  ,\nonumber\\
\partial_{x}v(0,t)  &  =\partial_{x}v(1,t)=0,\text{ \ \ \ }t\in\left[
0,T\right)  , \label{example2backward}%
\end{align}
where $\varphi$ and $\psi$ satisfy the hypotheses of the preceding section.

As is well known, the solution to (\ref{example1forward}) may be written as%
\begin{equation}
u_{\varphi}(x,t)=\sum_{n\in\mathbb{Z}}a_{n}\cos\left(  \pi nx\right)
\exp\left[  -\frac{\pi^{2}n^{2}}{2}t\right]  \label{series1}%
\end{equation}
with%
\[
a_{n}=\int_{0}^{1}dx\varphi(x)\cos\left(  \pi nx\right)
\]
for every $n\in\mathbb{Z}$, while the solution to (\ref{example2backward})
reads%
\begin{equation}
v_{\psi}(x,t)=\sum_{n\in\mathbb{Z}}b_{n}\cos\left(  \pi nx\right)  \exp\left[
-\frac{\pi^{2}n^{2}}{2}\left(  T-t\right)  \right]  \label{series2}%
\end{equation}
with%
\[
b_{n}=\int_{0}^{1}dx\psi(x)\cos\left(  \pi nx\right)  ,
\]
both series (\ref{series1}) and (\ref{series2}) being absolutely and uniformly
convergent. Furthermore (\ref{evolutionsystem1}) and (\ref{evolutionsystem2})
hold with the kernels%
\begin{align}
&  g(x,t;y,s)\nonumber\\
&  =g^{\ast}\left(  x,s;y,t\right)  =\sum_{n\in\mathbb{Z}}\cos\left(  \pi
nx\right)  \cos\left(  \pi ny\right)  \exp\left[  -\frac{\pi^{2}n^{2}}%
{2}(t-s)\right]  , \label{onedimgreen}%
\end{align}
and moreover all of the results of the preceding section are valid providing
one chooses a joint probability distribution of the form (\ref{specialclass})
satisfying (\ref{normalizationcondition}). Equations (\ref{itoprocess1}) and
(\ref{itoprocess2}) then read%
\begin{equation}
Z_{t}=Z_{0}+\int_{0}^{t}d\tau\partial_{x}\ln v_{\psi}(Z_{\tau},\tau
)+W_{t}^{\ast} \label{itoprocess5}%
\end{equation}
and%
\begin{equation}
Z_{t}=Z_{T}+\int_{t}^{T}d\tau\partial_{x}\ln u_{\varphi}(Z_{\tau},\tau)+W_{t},
\label{itoprocess6}%
\end{equation}
respectively.

Our goal now is to choose as simple a $\varphi$ and $\psi$ as possible in
order to unveil further properties of the processes thus constructed, keeping
in mind that we ought to disregard the trivial case%
\[
\varphi(x)=\psi(x)=1,
\]
which does imply (\ref{normalizationcondition}) but gives%
\[
u_{\varphi}(x,t)=v_{\psi}(x,t)=1
\]
for all $(x,t)\in\left[  0,1\right]  \times\left[  0,T\right]  $ and thereby
$b^{\ast}=b=0$ identically according to (\ref{drift4}) and (\ref{drift5}). Let
us choose instead%
\begin{equation}
\varphi(x)=1+\frac{1}{2}\cos\left(  \pi x\right)  \label{initialdatum}%
\end{equation}
and%
\begin{equation}
\psi(x)=1. \label{finaldatum}%
\end{equation}
Then the normalization condition (\ref{normalizationcondition}) holds and we
have%
\begin{equation}
u_{\varphi}(x,t)=1+\frac{1}{2}\cos\left(  \pi x\right)  \exp\left[  -\frac
{\pi^{2}}{2}t\right]  \label{forwardsolution2}%
\end{equation}
and%
\begin{equation}
v_{\psi}(x,t)=1 \label{backwardsolution2}%
\end{equation}
for all $(x,t)\in\left[  0,1\right]  \times\left[  0,T\right]  $. We have the
following result:

\bigskip

\textbf{Proposition 4.} \textit{Let us consider (\ref{example1forward}) and
(\ref{example2backward}) with the data (\ref{initialdatum}) and
(\ref{finaldatum}), respectively. Then the following statements hold:}

\textit{(a) We have}%
\[
\mathbb{P}_{\mu}\left(  Z_{t}\in E\right)  =\left\vert E\right\vert +\frac
{1}{2}\exp\left[  -\frac{\pi^{2}}{2}t\right]  \int_{E}dx\cos\left(  \pi
x\right)
\]
\textit{for every Borel subset }$E$\textit{ }$\subseteq\left[  0,1\right]
$\textit{ of Lebesgue measure }$\left\vert E\right\vert $\textit{ and every
}$t\in\left[  0,T\right]  $\textit{.}

\textit{(b) When considered as a forward Markov diffusion the process
}$Z_{\tau\in\left[  0,T\right]  }$ \textit{is a Wiener process with zero drift
which is instantaneously reflected at }$x=0$\textit{ and }$x=1$.\textit{
Moreover, the Lebesgue measure on }$\left[  0,1\right]  $ \textit{is}
\textit{an invariant measure for }$Z_{\tau\in\left[  0,T\right]  }$.

\textit{(c) When considered as a backward Markov diffusion the process
}$Z_{\tau\in\left[  0,T\right]  }$ \textit{satisfies}%
\[
Z_{t}=Z_{T}-\pi%
{\displaystyle\int\limits_{t}^{T}}
d\tau\frac{\sin\left(  \pi Z_{\tau}\right)  \exp\left[  -\frac{\pi^{2}}{2}%
\tau\right]  }{2+\cos\left(  \pi Z_{\tau}\right)  \exp\left[  -\frac{\pi^{2}%
}{2}\tau\right]  }+W_{t}%
\]
\textit{for every} $t\in\left[  0,T\right]  $\textit{, which corresponds to
the backward drift}%
\begin{equation}
b(x,t)=\frac{\pi\sin\left(  \pi x\right)  \exp\left[  -\frac{\pi^{2}}%
{2}t\right]  }{2+\cos\left(  \pi x\right)  \exp\left[  -\frac{\pi^{2}}%
{2}t\right]  } \label{onedimbackdrift}%
\end{equation}
\textit{defined} \textit{for all} $(x,t)\in\left[  0,1\right]  \times\left[
0,T\right]  $\textit{.}

\bigskip

\textbf{Proof. }Statement (a) follows from (\ref{bernsteindensity}),
(\ref{forwardsolution2}) and (\ref{backwardsolution2}).

Let us now consider $Z_{\tau\in\left[  0,T\right]  }$ as a forward Markov
diffusion. Then we have%
\[
b^{\ast}(x,t)=0
\]
for all $(x,t)\in\left[  0,1\right]  \times\left[  0,T\right]  $ according to
(\ref{drift4}), so that (\ref{itoprocess5}) reduces to%
\[
Z_{t}=Z_{0}+W_{t}^{\ast}%
\]
which defines a Wiener process with zero drift. As for the instantaneous
reflection at the boundaries of the domain, let us condition the process by
setting $Z_{0}=x$ for an arbitrary $x\in\left(  0,1\right)  $. Since $Z_{t}%
\in\left[  0,1\right]  $ we then have%
\[
Z_{t}=\left\vert x+W_{t}^{\ast}\right\vert
\]
and by virtue of the identity%
\begin{align*}
&  \sum_{n\in\mathbb{Z}}\cos\left(  \pi nx\right)  \cos\left(  \pi ny\right)
\exp\left[  -\frac{\pi^{2}n^{2}}{2}t\right] \\
&  =\left(  2\pi t\right)  ^{-\frac{1}{2}}\sum_{n\in\mathbb{Z}}\left(
\exp\left[  -\frac{\left\vert x+y+2n\right\vert ^{2}}{2t}\right]  +\exp\left[
-\frac{\left\vert x-y-2n\right\vert ^{2}}{2t}\right]  \right)
\end{align*}
valid for all $x,y\in\left(  0,1\right)  $ and every $t\in\left(  0,T\right)
$ (see, for instance, Appendix 1 to Chapter 6 in \cite{borodinsalminen}), we
have%
\begin{align*}
&  \mathbb{P}_{\mu}^{x}\left(  Z_{t}\in E\right) \\
&  =\int_{E}dyg^{\ast}(x,0;y,t)\\
&  =\left(  2\pi t\right)  ^{-\frac{1}{2}}\int_{E}dy\sum_{n\in\mathbb{Z}%
}\left(  \exp\left[  -\frac{\left\vert x+y+2n\right\vert ^{2}}{2t}\right]
+\exp\left[  -\frac{\left\vert x-y-2n\right\vert ^{2}}{2t}\right]  \right)
\end{align*}
for each Borel subset $E\subseteq\left[  0,1\right]  $ and all $(x,t)\in
\left(  0,1\right)  \times\left(  0,T\right)  $. Therefore, we may indeed
identify $Z_{\tau\in\left[  0,T\right]  }$ with the doubly reflected Brownian
motion as defined for instance in Section 2.8 C of \cite{karatzasshreve}.
Finally let us consider the function (\ref{forwardtransition}), which in this
case takes the form%
\[
M^{\ast}\left(  x,s;E,t\right)  =\int_{E}dyg^{\ast}(x,s;y,t)
\]
where $g^{\ast}$ is given by (\ref{onedimgreen}). Since we have%
\[
\int_{0}^{1}dxM^{\ast}\left(  x,s;E,t\right)  =\left\vert E\right\vert
\]
for each Borel subset $E\subseteq\left[  0,1\right]  $ and all $s,t\in\left[
0,T\right]  $ satisfying $t>s$, the Lebesgue measure is indeed an invariant
measure for $Z_{\tau\in\left[  0,T\right]  }$ on $\left[  0,1\right]  $.
Consequently Statement (b)\ holds.

Statement (c) follows directly from (\ref{drift5}), (\ref{itoprocess6}) and
(\ref{forwardsolution2}). \ \ $\blacksquare$

\bigskip

\textsc{Remark.} The Bernstein diffusion $Z_{\tau\in\left[  0,T\right]  }$ of
the preceding proposition is not reversible in the sense of \cite{dobsukfritz}
or \cite{kolmogorov} since its probability density is time-dependent according
to Statement (a). This is a simple consequence of the third remark following
the proof of Theorem 2.

\bigskip

\textsc{Example 2.} We now consider the initial-boundary value problem%
\begin{align}
\partial_{t}u(x,t)  &  =\frac{1}{2}\triangle_{x}u(x,t),\text{ \ \ }%
(x,t)\in\mathbb{D}\times\left(  0,T\right]  ,\text{\ }\nonumber\\
u(x,0)  &  =\varphi(x),\text{ \ \ }x\in\mathbb{D},\nonumber\\
\frac{\partial u(x,t)}{\partial n(x)}  &  =0,\text{ \ \ }(x,t)\in
\partial\mathbb{D\times}\left(  0,T\right]  \label{example3forward}%
\end{align}
and the corresponding final-boundary value problem for the adjoint equation%
\begin{align}
-\partial_{t}v(x,t)  &  =\frac{1}{2}\triangle_{x}v(x,t),\text{ \ \ }%
(x,t)\in\mathbb{D}\times\left[  0,T\right)  ,\nonumber\\
v(x,T)  &  =\psi(x),\text{ \ \ }x\in\mathbb{D},\nonumber\\
\frac{\partial v(x,t)}{\partial n(x)}  &  =0,\text{ \ \ \ }(x,t)\in
\partial\mathbb{D\times}\left[  0,T\right)  , \label{example4backward}%
\end{align}
where%
\[
\mathbb{D=}\left\{  x\in\mathbb{R}^{2}:\left\vert x\right\vert <1\right\}
\]
is the two-dimensional open unit disk centered at the origin, and where
$\varphi$ and $\psi$ satisfy the hypotheses of the preceding section. Let us
restrict ourselves to radially symmetric solutions to (\ref{example3forward})
and (\ref{example4backward}); switching to polar coordinates and abusing the
notation a bit we then rewrite (\ref{example3forward}) as%
\begin{align}
\partial_{t}u(r,t)  &  =\frac{1}{2}\left(  \partial_{rr}+r^{-1}\partial
_{r}\right)  u(r,t),\text{ \ \ }(r,t)\in\left(  0,1\right]  \times\left(
0,T\right]  ,\text{\ }\nonumber\\
u(r,0)  &  =\varphi(r),\text{ \ \ }r\in\left[  0,1\right]  ,\nonumber\\
\partial_{r}u(1,t)  &  =0,\text{ \ \ }t\in\left(  0,T\right]  ,
\label{polarcoordinates1}%
\end{align}
whose solution reads%
\begin{equation}
u_{\varphi}(r,t)=\sum_{n=1}^{+\infty}a_{n,0}J_{0}\left(  \sqrt{\mu_{n,0}%
}r\right)  \exp\left[  -\frac{\mu_{n,0}}{2}t\right]  \label{series3}%
\end{equation}
where $J_{0}$ stands for the Bessel function of zeroth order of the first kind
(see, for instance, \cite{watson} or \cite{weinberger}). In (\ref{series3}) we
have%
\[
a_{n,0}=2J_{0}^{-2}\left(  \sqrt{\mu_{n,0}}\right)  \int_{0}^{1}%
rdr\varphi\left(  r\right)  J_{0}\left(  \sqrt{\mu_{n,0}}r\right)
\]
where the $\mu_{n,0}$'s are ordered in such a way that%
\begin{equation}
0=\mu_{1,0}<\mu_{2,0}<\mu_{3,0}<..... \label{arrangement}%
\end{equation}
and satisfy%
\begin{equation}
J_{0}^{\prime}\left(  \sqrt{\mu_{n,0}}\right)  =0, \label{neumanncondition}%
\end{equation}
the preceding relation being Neumann's boundary condition in this case. In a
similar way we have%
\begin{align}
-\partial_{t}v(r,t)  &  =\frac{1}{2}\left(  \partial_{rr}+r^{-1}\partial
_{r}\right)  v(r,t),\text{ \ \ }(r,t)\in\left(  0,1\right]  \times\left[
0,T\right)  ,\nonumber\\
v(r,T)  &  =\psi(r),\text{ \ \ }r\in\left[  0,1\right]  ,\nonumber\\
\partial_{r}v(1,t)  &  =0,\text{ \ \ }t\in\left[  0,T\right)  ,
\label{polarcoordinates2}%
\end{align}
for (\ref{example4backward}), whose solution is%
\begin{equation}
v_{\psi}(r,t)=\sum_{n=1}^{+\infty}b_{n,0}J_{0}\left(  \sqrt{\mu_{n,0}%
}r\right)  \exp\left[  -\frac{\mu_{n,0}}{2}\left(  T-t\right)  \right]
\label{series4}%
\end{equation}
where%
\[
b_{n,0}=2J_{0}^{-2}\left(  \sqrt{\mu_{n,0}}\right)  \int_{0}^{1}rdr\psi\left(
r\right)  J_{0}\left(  \sqrt{\mu_{n,0}}r\right)  ,
\]
both series (\ref{series3}) and (\ref{series4})\ being again absolutely and
uniformly convergent. Consequently we may write%
\[
u_{\varphi}(r,t)=\int_{0}^{1}r^{\prime}dr^{\prime}g\left(  r,t;r^{\prime
},0\right)  \varphi\left(  r^{\prime}\right)
\]
and%
\[
v_{\psi}(r,t)=\int_{0}^{1}r^{\prime}dr^{\prime}g^{\ast}\left(  r,t;r^{\prime
},T\right)  \psi\left(  r^{\prime}\right)
\]
where the radial Green functions are%
\begin{align*}
&  g\left(  r,t;r^{\prime},s\right) \\
&  =g^{\ast}\left(  r,s;r^{\prime},t\right) \\
&  =2\sum_{n=1}^{+\infty}J_{0}^{-2}\left(  \sqrt{\mu_{n,0}}\right)
J_{0}\left(  \sqrt{\mu_{n,0}}r\right)  J_{0}\left(  \sqrt{\mu_{n,0}}r^{\prime
}\right)  \exp\left[  -\frac{\mu_{n,0}}{2}\left(  t-s\right)  \right]  ,
\end{align*}
so that all of the results of the preceding section remain valid in this case
provided one chooses again a probability distribution of the form
(\ref{specialclass}) satisfying (\ref{normalizationcondition}).

Let us take for instance%
\begin{equation}
\varphi(r)=\frac{1}{\pi}\left(  1+J_{0}\left(  \sqrt{\mu_{2,0}}r\right)
\right)  \label{initialdatum2}%
\end{equation}
and%
\begin{equation}
\psi\left(  r\right)  =1 \label{finaldatum2}%
\end{equation}
for every $r\in\left[  0,1\right]  $, where $\sqrt{\mu_{2,0}}$ is the first
positive root of $J_{0}^{\prime}$ according to (\ref{arrangement}) and
(\ref{neumanncondition}). It then follows from the standard properties of
$J_{0}$ that $\varphi$ satisfies all the required hypotheses of our theory,
including positivity and (\ref{normalizationcondition}); furthermore, the
orthogonality properties of $J_{0}$ with respect ot the $\mu_{n,0}$'s imply
that the corresponding solutions are%
\begin{equation}
u_{\varphi}(r,t)=\frac{1}{\pi}\left(  1+J_{0}\left(  \sqrt{\mu_{2,0}}r\right)
\exp\left[  -\frac{\mu_{2,0}}{2}t\right]  \right)  \label{forwardsolution3}%
\end{equation}
and%
\begin{equation}
v_{\psi}(r,t)=1 \label{backwardsolution3}%
\end{equation}
for all $(r,t)\in\left[  0,1\right]  \times\left[  0,T\right]  $. Let us
denote again by $Z_{\tau\in\left[  0,T\right]  }$ the Bernstein diffusion
associated with (\ref{polarcoordinates1}), (\ref{polarcoordinates2}) and
(\ref{forwardsolution3}), (\ref{backwardsolution3}), respectively, and by
$J_{1}=-J_{0}^{\prime}$ the Bessel function of first order. We have the
following result, similar to that of Proposition 4:

\bigskip

\textbf{Proposition 5.}\textit{ Let us consider (\ref{polarcoordinates1}) and
(\ref{polarcoordinates2}) with the data (\ref{initialdatum2}) and
(\ref{finaldatum2}), respectively. Then the following statements hold:}

\textit{(a) We have}%
\[
\mathbb{P}_{\mu}\left(  Z_{t}\in E\right)  =\frac{1}{\pi}\left(  \left\vert
E\right\vert +\exp\left[  -\frac{\mu_{2,0}}{2}t\right]  \int_{E}dxJ_{0}\left(
\sqrt{\mu_{2,0}}\left\vert x\right\vert \right)  \right)
\]
\textit{for every Borel subset }$E$\textit{ }$\subseteq\mathbb{D}$\textit{ of
Lebesgue measure }$\left\vert E\right\vert $\textit{ and every }$t\in\left[
0,T\right]  $\textit{.}

\textit{(b) When considered as a forward Markov diffusion the process
}$Z_{\tau\in\left[  0,T\right]  }$ \textit{is a Wiener process with zero
drift} \textit{and} \textit{the Lebesgue measure on }$\mathbb{D}$ \textit{is}
\textit{an invariant measure for }$Z_{\tau\in\left[  0,T\right]  }$.

\textit{(c) When considered as a backward Markov diffusion the process
}$Z_{\tau\in\left[  0,T\right]  }$ \textit{satisfies}%
\[
Z_{t}=Z_{T}-%
{\displaystyle\int\limits_{t}^{T}}
d\tau b\left(  Z_{\tau},\tau\right)  +W_{t}%
\]
\textit{for every} $t\in\left[  0,T\right]  $\textit{, where the backward
drift is given by}%
\begin{equation}
b(x,t)=\frac{\sqrt{\mu_{2,0}}J_{1}\left(  \sqrt{\mu_{2,0}}\left\vert
x\right\vert \right)  \exp\left[  -\frac{\mu_{2,0}}{2}t\right]  }%
{1+J_{0}\left(  \sqrt{\mu_{2,0}}\left\vert x\right\vert \right)  \exp\left[
-\frac{\mu_{2,0}}{2}t\right]  }\times\frac{x}{\left\vert x\right\vert }
\label{backdrift1}%
\end{equation}
\textit{for all} $(x,t)\in\mathbb{D}\times\left[  0,T\right]  $\textit{ with
}$x\neq0$\textit{, and by}%
\begin{equation}
b(0,t)=\lim_{\left\vert x\right\vert \rightarrow0}b(x,t)=0 \label{backdrift2}%
\end{equation}
\textit{for every} $t\in\left[  0,T\right]  $.

\bigskip

\textbf{Proof. }Statement (a) follows from (\ref{bernsteindensity}), this time
with (\ref{forwardsolution3}) and (\ref{backwardsolution3}).

Let us now consider $Z_{\tau\in\left[  0,T\right]  }$ as a forward Markov
diffusion. Then we have%
\[
Z_{t}=Z_{0}+W_{t}^{\ast}%
\]
since%
\[
b^{\ast}(x,t)=0
\]
for all $(x,t)\in\mathbb{D}\times\left[  0,T\right]  $ according to
(\ref{drift4}) and (\ref{backwardsolution3}), which indeed defines a Wiener
process with zero drift, and the fact that the Lebesgue measure on
$\mathbb{D}$ is an invariant measure for $Z_{\tau\in\left[  0,T\right]  }$
follows from an argument similar to that given in Proposition 4. Consequently
Statement (b)\ holds.

Statement (c) is a simple consequence of (\ref{drift5}), (\ref{itoprocess2})
and (\ref{forwardsolution3}), with (\ref{backdrift2}) following from the fact
that $J_{0}(0)=1$ and%
\[
J_{1}\left(  \sqrt{\mu_{2,0}}\left\vert x\right\vert \right)  \sim\frac{1}%
{2}\sqrt{\mu_{2,0}}\left\vert x\right\vert
\]
as $\left\vert x\right\vert \rightarrow0$. \ \ $\blacksquare$

\bigskip

\textsc{Remarks.} (1) Relations (\ref{backdrift1}) and (\ref{backdrift2})
reveal two interesting features regarding the drift of the backward diffusion
$Z_{\tau\in\left[  0,T\right]  }$. On the one hand, $b(x,t)$ points to the
radial direction and its norm $\left\vert b(x,t)\right\vert $ is rotationally
invariant, which is no surprise since (\ref{forwardsolution3}) is radially
symmetric to start with. On the other hand, $\left\vert b(x,t)\right\vert $
vanishes at the center of the disk and on its boundary since $J_{1}\left(
\sqrt{\mu_{2,0}}\right)  =0$, and thereby reaches a maximal value at some
$\left\vert x^{\ast}\right\vert =r^{\ast}\in\left(  0,1\right)  $ depending on
$t$ since $J_{1}>0$ on the interval $\left(  0,\sqrt{\mu_{2,0}}\right)  $. We
note that such a feature is already present in the one-dimensional backward
drift (\ref{onedimbackdrift}).

(2) The differential operator on the right-hand side of
(\ref{polarcoordinates1}) and (\ref{polarcoordinates2}) is, of course,
formally the generator of a Bessel process of order zero. Along the same line
we can also consider radially symmetric solutions to Problems
(\ref{example3forward}) and (\ref{example4backward}) in the $d$-dimensional
ball centered at the origin of $\mathbb{R}^{d}$, where the radial part of half
of the Laplacian identifies with the generator of a Bessel process of order
$\nu:=\frac{d}{2}-1$. In such cases we can naturally expect the existence of
nontrivial connections between Bernstein diffusions and Bessel processes of
order $\nu$. Moreover, because of (\ref{boundarycondition1}) and
(\ref{boundarycondition2}) it is tempting to conjecture that all Bernstein
diffusions constructed in this article are in fact reflected diffusions, as is
the case in Statement (b) of Proposition 4. We defer a detailed analysis of
these questions to a separate publication, including one of the boundary local
time%
\[
L_{Z}(t):=\int_{0}^{t}ds\mathbb{I}_{\partial\mathbb{D}}(Z_{s})=\lim
_{\varepsilon\rightarrow0_{+}}\frac{1}{\varepsilon}\int_{0}^{t}ds\mathbb{I}%
_{\mathbb{D}_{\varepsilon}}(Z_{s})
\]
in light of the results of \cite{lionsznitman}, \cite{satoueno} and
\cite{tanaka}, where $\mathbb{I}_{\partial\mathbb{D}}$ is the indicator
function of the boundary $\partial\mathbb{D}$ and $\mathbb{I}_{\mathbb{D}%
_{\varepsilon}}$that of the annulus
\[
\mathbb{D}_{\varepsilon}:=\left\{  x\in\mathbb{R}^{d}:0<1-\varepsilon
\leq\left\vert x\right\vert \leq1\right\}  .
\]

(3) We conclude this article by observing that problems of the form
(\ref{parabolicproblem}) and (\ref{adjointproblem}) can play an important
r\^{o}le in physics and in the natural sciences in general, as they encompass
equations which are quite relevant to the mathematical analysis of various
phenomena ranging from the propagation of heat to the space-time evolution of
gene densities in population dynamics or population genetics, to name only a
few (see, for instance, \cite{chueshovvuiller} and the references therein).
The probabilistic interpretation we have developed in this article may,
therefore, also shed new light on such phenomena, a remark that applies in
particular to problems (\ref{example1forward})-(\ref{example2backward}) and
(\ref{example3forward})-(\ref{example4backward}), respectively. Furthermore,
with $i=\sqrt{-1}$ we also see that the formal change of variables $t\mapsto
it$ reduces the equations of Section 3 to the basic Schr\"{o}dinger equations
of quantum physics. For instance, with a slight abuse of notation again, the
first equation in (\ref{example3forward}) becomes%
\[
i\partial_{t}u(x,t)=-\frac{1}{2}\triangle_{x}u(x,t),\text{ \ \ }%
(x,t)\in\mathbb{D}\times\left(  0,T\right]  ,\text{\ }%
\]
while (\ref{example4backward}) transforms into%
\[
-i\partial_{t}v(x,t)=-\frac{1}{2}\triangle_{x}v(x,t),\text{ \ \ }%
(x,t)\in\mathbb{D}\times\left[  0,T\right)
\]
with $v=\bar{u}$ in case of uniqueness. Consequently, the usual
quantum-mechanical interpretation of the quantity%

\[
\int_{E}dx\left\vert u(x,t)\right\vert ^{2}%
\]
as the probability of finding a particle in the region $E$ at time $t$ has a
rigorous parabolic counterpart in the form of Statement (c) in Theorem 2. 

\textbf{Acknowledgments.} The research of both authors was supported by the
FCT\textit{ }of the Portuguese government under grant PTDC/MAT/69635/2006, and
by the Mathematical Physics Group\textit{ }of the University of Lisbon under
grant ISFL/1/208. The first author is also indebted to Madalina Deaconu and
Elton Hsu for stimulating discussions and correspondence on the theme of
reflected diffusions. Last but not least, he wishes to thank the Complexo
Interdisciplinar da Universidade de Lisboa and the ETH-Forschungsinstitut
f\"{u}r Mathematik in Zurich\textit{ }where parts of this work were completed
for their financial support and warm hospitality.

\section{Appendix. A Variational Construction of Weak Solutions in $L^{2}%
(D)$.}

The solutions $u_{\varphi}$ and $v_{\psi}$ we used throughout this article are
generated by two evolution systems $U_{A}(t,s)_{0\leq s\leq t\leq T}$ and
$U_{A}^{\ast}(t,s)_{0\leq s\leq t\leq T}$ on $L^{2}(D)$. We show here how to
construct these evolution systems by applying the standard methods of
\cite{tanabe}, under the following hypotheses regarding the coefficients $k$,
$l$ and $V$ in (\ref{parabolicproblem}) and (\ref{adjointproblem}):

\bigskip

(K$^{\prime}$) The function $k:D\times\left[  0,T\right]  \mapsto
\mathbb{R}^{d^{2}}$ is matrix-valued and for every $i,j\in\left\{
1,...,d\right\}  $ we have $k_{i,j}=k_{j,i}\in L^{\infty}(D\times\left(
0,T\right)  \mathbb{)}$;\ moreover, there exists a finite constant
$\underline{k}>0$ such that the inequality%
\begin{equation}
\left(  k(x,t)q,q\right)  _{\mathbb{R}^{d}}\geq\underline{k}\left\vert
q\right\vert ^{2} \label{ellipticitycondition}%
\end{equation}
holds uniformly in $(x,t)\in D\times\left[  0,T\right]  $ for all
$q\in\mathbb{R}^{d}$. Finally, there exist finite constants $c_{\ast}>0$,
$\beta\in\left(  \frac{1}{2},1\right]  $ such that the H\"{o}lder continuity
estimate%
\[
\max_{i,j\in\left\{  1,...,d\right\}  }\left\vert k_{i,j}(x,t)-k_{i,j}%
(x,s)\right\vert \leq c_{\ast}\left\vert t-s\right\vert ^{\beta}%
\]
is valid for every $x\in D$ and all $s,t\in\left[  0,T\right]  $.

\bigskip

As for the lower-order differential operators we assume that the following
hypotheses are valid, where we assume without restricting the generality that
the constants $c_{\ast}$ and $\beta$ are the same as in Hypothesis
(K$^{\prime}$):

\bigskip

(L$^{\prime}$) Each component of the vector-field $l:D\times\left[
0,T\right]  \mapsto\mathbb{R}^{d}$ satisfies $l_{i}\in L^{\infty}%
(D\times\left(  0,T\right)  \mathbb{)}$. Moreover, the H\"{o}lder continuity
estimate%
\[
\max_{i\in\left\{  1,...,d\right\}  }\left\vert l_{i}(x,t)-l_{i}%
(x,s)\right\vert \leq c_{\ast}\left\vert t-s\right\vert ^{\beta}%
\]
holds for every $x\in D$ and all $s,t\in\left[  0,T\right]  $.

\bigskip

(V$^{\prime}$) The function $V:D\times(0,T)\mapsto\mathbb{R}$ is such that
$V\in L^{\infty}(D\times(0,T))$ and satisfies%
\[
\left\vert V(x,t)-V(x,s)\right\vert \leq c_{\ast}\left\vert t-s\right\vert
^{\beta}%
\]
for every $x\in D$ and all $s,t\in\left[  0,T\right]  $.

\bigskip

Moreover, both the initial condition $\varphi$ and the final condition $\psi$
are real-valued and the following hypothesis holds:

\bigskip

(IF$^{\prime}$) We have $\varphi,\psi\in L^{2}(D)$.

\bigskip

\textsc{Remark.} In the variational theory we are reviewing here we observe
that the H\"{o}lder continuity requirement relative to the time variable in
Hypotheses (K$^{\prime}$), (L$^{\prime}$) and (V$^{\prime}$) is stronger than
that of Hypotheses (K), (L) and (V), since $\beta\in\left(  \frac{1}%
{2},1\right]  $ whereas $\frac{\alpha}{2}\in\left(  0,\frac{1}{2}\right)  $.
However, it is easy to show by uniqueness arguments that the evolution
operators $U_{A}(t,s)_{0\leq s\leq t\leq T}$ and $U_{A}^{\ast}(t,s)_{0\leq
s\leq t\leq T}$ introduced in Section 2 are identical to those constructed
below. The reason why $\beta\in\left(  \frac{1}{2},1\right]  $ is required
here is intimately tied up with the variational structure of the problem, and
is thoroughly discussed in \cite{tanabe}.

\bigskip

Under the preceding three conditions, it is easily verified that the quadratic
form $a:\left[  0,T\right]  \times H^{1}(D\mathbb{)\times}H^{1}%
(D\mathbb{)\mapsto C}$ defined by%
\begin{align*}
a\left(  t,f,h\right)   &  :=\frac{1}{2}\int_{D}dx\left(  k(x,t)\nabla
_{x}f(x),\nabla_{x}h(x)\right)  _{\mathbb{C}^{d}}\\
&  +\int_{D}dx\left(  l(x,t),\nabla_{x}f(x)\right)  _{\mathbb{C}^{d}}%
\overline{h}(x)\\
&  +\int_{D}dxV(x,t)f(x)\overline{h}(x)
\end{align*}
satisfies the estimates%
\begin{align}
\left\vert a\left(  t,f,h\right)  \right\vert  &  \leq c\left\Vert
f\right\Vert _{1,2}\left\Vert h\right\Vert _{1,2},\label{estimate1}\\
\operatorname{Re}a\left(  t,f,f\right)   &  \geq\underline{k}\left\Vert
f\right\Vert _{1,2}^{2}-c\left\Vert f\right\Vert _{2}^{2},\label{estimate2}\\
\left\vert a\left(  t,f,h\right)  -a\left(  s,f,h\right)  \right\vert  &  \leq
c\left\vert t-s\right\vert ^{\beta}\left\Vert f\right\Vert _{1,2}\left\Vert
h\right\Vert _{1,2} \label{estimate3}%
\end{align}
for all $s,t\in\left[  0,T\right]  $ and all $f,h\in H^{1}(D\mathbb{)}$, where
$\left\Vert .\right\Vert _{2}$ and $\left\Vert .\right\Vert _{1,2}$ stand for
the usual norms in $L^{2}(D\mathbb{)}$ and $H^{1}(D\mathbb{)}$, respectively,
and where $\left(  .,.\right)  _{\mathbb{C}^{d}}$ denotes the Hermitian inner
product in $\mathbb{C}^{d}$. Consequently, the formal elliptic operator%
\[
A(t):=-\frac{1}{2}\operatorname{div}\left(  k(.,t)\nabla\right)  +\left(
l(.,t),\nabla\right)  _{\mathbb{C}^{d}}+V(.,t)
\]
corresponding to the right-hand side of (\ref{parabolicproblem}) can be
realized as a regularly accretive operator defined on some time-dependent and
dense domain $\mathcal{D(}A(t))\subset L^{2}(D\mathbb{)}$, and as such
generates an evolution system $U_{A}(t,s)_{0\leq s\leq t\leq T}$ in
$L^{2}(D\mathbb{)}$ given by%
\begin{equation}
U_{A}(t,s)f(x)=\left\{
\begin{array}
[c]{c}%
f(x)\text{ \ \ \ \ \ \ \ \ \ \ \ \ \ \ \ \ \ \ \ \ \ \ \ \ if }t=s,\\
\int_{D}dyg_{A}(x,t;y,s)f(y)\text{ \ \ if }t>s
\end{array}
\right.  \label{evolution}%
\end{equation}
for every $f\in L^{2}(D\mathbb{)}$, where $g_{A}$ denotes the parabolic Green
function associated with (\ref{parabolicproblem}). Indeed all these assertions
follow directly from estimates (\ref{estimate1})-(\ref{estimate3}) and the
general theory developed in Section 5.4 of \cite{tanabe}, together with
Schwartz's kernel theorem which guarantees the existence of $g_{A}$ (see
\cite{schwartz} for a summary of the many possible applications of that theorem).

In a similar way, the Hermitian conjugate form%
\[
a^{\ast}\left(  t,f,h\right)  :=\overline{a\left(  t,h,f\right)  }%
\]
is associated with the linear operator $A^{\ast}(t)$ adjoint to $A(t)$, which
in turn generates the adjoint evolution system%
\begin{equation}
U_{A}^{\ast}(t,s)f(x)=\left\{
\begin{array}
[c]{c}%
f(x)\text{ \ \ \ \ \ \ \ \ \ \ \ \ \ \ \ \ \ \ \ \ \ \ \ \ if }t=s,\\
\int_{D}dyg_{A}^{\ast}(x,s;y,t)f(y)\text{ \ \ if }t>s,
\end{array}
\right.  \label{adjointevolution}%
\end{equation}
where $G_{A}^{\ast}$ is the parabolic Green function associated with
(\ref{adjointproblem}) that satisfies the relation%
\[
g_{A}^{\ast}(x,s;y,t)=g_{A}(y,t;x,s)
\]
for all $s,t\in\left[  0,T\right]  $ with $t>s$.

The important features of (\ref{evolution}) and (\ref{adjointevolution}) are
that they provide the real-valued functions defined by%
\begin{equation}
u_{\varphi}(x,t):=U_{A}(t,0)\varphi(x)=\int_{D}dyg_{A}(x,t;y,0)\varphi
(y),\text{ \ \ }t\in\left(  0,T\right]  \label{weakforwardsolution}%
\end{equation}
and%
\begin{equation}
v_{\psi}(x,t):=U_{A}^{\ast}(T,t)\psi(x)=\int_{D}dyg_{A}^{\ast}(x,t;y,T)\psi
(y),\text{ \ \ }t\in\left[  0,T\right)  , \label{weakbackwardsolution}%
\end{equation}
which satisfy$\ \ \ \ \ \ \ \ \ \ \ $%
\[
\left(  \frac{\partial}{\partial t}u_{\varphi}(.,t),h\right)  _{2}%
+a(t,u_{\varphi}(.,t),h)=0,\text{ \ \ }t\in\left(  0,T\right]
\]
and%
\[
-\left(  \frac{\partial}{\partial t}v_{\psi}(.,t),h\right)  _{2}+a^{\ast
}(t,v_{\psi}(.,t),h)=0,\text{ \ \ }t\in\left[  0,T\right)
\]
for every $h\in H^{1}(D\mathbb{)}$, respectively, where $(.,.)_{2}$ stands for
the usual inner product in $L^{2}(D\mathbb{)}$. Moreover we have $u_{\varphi
},v_{\psi}\in L^{2}(D\times(0,T))$, so that\ (\ref{weakforwardsolution}) and
(\ref{weakbackwardsolution}) provide \textit{weak} solutions to
(\ref{parabolicproblem}) and (\ref{adjointproblem}), respectively (see, for
instance, Section 5.5 in \cite{tanabe}).

These solutions are those which ultimately possess the properties listed in
Lemma 1 of Section 2, according to the above remark regarding the H\"{o}lder
regularity in time.


\begin{thebibliography}{99}                                                                                               %


\bibitem {adamsfournier}\textsc{Adams, R. A., Fournier J. J. F.,}
\textit{Sobolev Spaces,} Pure Appl. Math. \textbf{140}, Academic Press, New
York (2003).

\bibitem {bernstein}\textsc{Bernstein, S.,} \textit{Sur les liaisons entre les
grandeurs al\'{e}atoires,} Verhandlungen des Internationalen
Mathematikerkongress \textbf{1} (1932) 288-309.

\bibitem {beurling}\textsc{Beurling, A.,} \textit{An Automorphism of Product
Measures,} Annals of Mathematics \textbf{72} (1960) 189-200.

\bibitem {borodinsalminen}\textsc{Borodin, A. N., Salminen, P.,}\textit{
Handbook of Brownian Motion-Facts and Formulae,} Probability and its
Applications Series, Birkh\"{a}user, Basel (2000).

\bibitem {carlen}\textsc{Carlen, E. A.,} \textit{Conservative Diffusions,}
Communications in Mathematical Physics \textbf{94} (1984)\ 293-315.

\bibitem {chueshovvuiller}\textsc{Chueshov, I. D., Vuillermot, P.-A.,}
\textit{Long-Time Behavior of Solutions to a Class of Quasilinear Parabolic
Equations with Random Coefficients}, Annales de l'Institut Henri
Poincar\'{e}-AN \textbf{15} (1998) 191-232.

\bibitem {cruzwuzambrin}\textsc{Cruzeiro, A. B.,} \textsc{Wu, L., Zambrini, J.
C.,} \textit{Bernstein Processes associated with Markov Processes, }in:
Stochastic Analysis and Mathematical Physics, (editor: R. Rebolledo)
Birkh\"{a}user, Basel (2000).

\bibitem {cruzeirozambrini}\textsc{Cruzeiro, A. B., Zambrini, J. C.,}
\textit{Malliavin Calculus and Euclidean Quantum Mechanics, I. Functional
Calculus, }Journal of Functional Analysis \textbf{96 }(1991) 62-95.

\bibitem {dobsukfritz}\textsc{Dobrushin, R. L., Sukhov, Y. M., Fritz, J.,}
\textit{A. N. Kolmogorov-The Founder of the Theory of Reversible Markov
Processes,} Russian Mathematical Surveys \textbf{43}, II (1988) 157-182.

\bibitem {dynkin}\textsc{Dynkin, E. B.,} \textit{Diffusions, Superdiffusions
and Partial Differential Equations,} American Mathematical Society Colloquium
Publications \textbf{50}, American Mathematical Society, Rhode Island (2002).

\bibitem {eideliva}\textsc{Eidelman, S. D., Ivasi\u{s}en, S. D.,}
\textit{Investigation of the Green Matrix for a Homogeneous Parabolic Boundary
Value Problem,} Transactions of the Moscow Mathematical Society \textbf{23}
(1970) 179-242.

\bibitem {eidelzhitara}\textsc{Eidelman, S. D., Zhitarashu, N. V.,}
\textit{Parabolic Boundary Value Problems,} Operator Theory, Advances and
Applications \textbf{101}, Birkh\"{a}user, Basel (1998).

\bibitem {freidlin}\textsc{Freidlin, M.,} \textit{Functional Integration and
Partial Differential Equations,} Annals of Mathematics Studies \textbf{109},
Princeton University Press (1985).

\bibitem {friedman}\textsc{Friedman, A.,} \textit{Partial Differential
Equations of Parabolic Type,} Prentice-Hall, Inc., Englewood Cliffs, New
Jersey (1964).

\bibitem {gihmanskohorod}\textsc{Gihman, I. I., Skorohod, A. V., }\textit{The
Theory of Stochastic Processes,} III, Classics in Mathematics Series,
Springer-Verlag, New York (2007).

\bibitem {gulicasteren}\textsc{Gulisashvili, A., van Casteren, J. A.,
}\textit{Non-Autonomous Kato Classes and Feynman-Kac Propagators,} World
Scientific, Singapore (2006).

\bibitem {hsu}\textsc{Hsu, P.,} \textit{Probabilistic Approach to the Neumann
Problem, }Communications on Pure and Applied Mathematics \textbf{38} (1985) 445-472.

\bibitem {ikedawatanabe}\textsc{Ikeda, N., Watanabe, S.,} \textit{Stochastic
Differential Equations and Diffusion Processes,} North-Holland Mathematical
Library \textbf{24}, North-Holland, Amsterdam (1989).

\bibitem {jamison}\textsc{Jamison, B.,} \textit{Reciprocal Processes,}
Zeitschrift f\"{u}r Wahrscheinlichkeitstheorie und Verwandte Gebiete
\textbf{30 }(1974) 65-86.

\bibitem {karatzasshreve}\textsc{Karatzas, I., Shreve, S. E.,}
\textit{Brownian Motion and Stochastic Calculus,} Graduate Texts in
Mathematics\textbf{ 113}, Springer-Verlag, New York (1987).

\bibitem {kolmogorov}\textsc{Kolmogorov, A. N.,} \textit{On the Reversibility
of the Statistical Laws of Nature,} in: Selected Works of A. N. Kolmogorov,
Mathematics and its Applications (Soviet Series) \textbf{26}, (editor: A. N.
Shiryayev) Kluwer, Boston (1992).

\bibitem {ladyurasolo}\textsc{Lady\u{z}enskaja, O. A., Solonnikov, V. A.,
Ural'ceva, N. N., }\textit{Linear and Quasilinear Equations of Parabolic
Type,} AMS Translations of Mathematical Monographs\textbf{ 23}, American
Mathematical Society, Providence, Rhode Island (1968).

\bibitem {lionsznitman}\textsc{Lions, P. L., Sznitman, A. S., }%
\textit{Stochastic Differential Equations with Reflecting Boundary
Conditions,} Communications on Pure and Applied Mathematics \textbf{37} (1984) 511-537.

\bibitem {nelson1}\textsc{Nelson, E.,} \textit{Dynamical Theories of Brownian
Motion,} Princeton University Press (1967).

\bibitem {privaultzambrin}\textsc{Privault, N., Zambrini, J. C.,}
\textit{Markovian Bridges and Reversible Diffusion Processes with Jumps,}
Annales de l'Institut Henri Poincar\'{e}-PR \textbf{40} (2004) 599-633.

\bibitem {satoueno}\textsc{Sato, K., Ueno, T.,}\textit{\ Multi-dimensional
Diffusion and the Markov Process on the Boundary,} Journal of Mathematics of
Kyoto University \textbf{4} (1964) 529-605.

\bibitem {schroedinger}\textsc{Schr\"{o}dinger, E.,} \textit{Sur la
th\'{e}orie relativiste de l'\'{e}lectron et l'interpr\'{e}tation de la
m\'{e}canique quantique,} Annales de l'Institut Henri Poincar\'{e} \textbf{2}
(1932) 269-310.

\bibitem {schwartz}\textsc{Schwartz, L.,} \textit{Th\'{e}orie des Noyaux,}
Proceedings of the International Congress of Mathematicians \textbf{1} (1950) 220-230.

\bibitem {solonnikov}\textsc{Solonnikov, V. A.,} \textit{On Boundary-Value
Problems for Linear Parabolic Systems of Differential Equations of General
Form,} Proceedings of the Steklov Institute of Mathematics \textbf{83} (1965).

\bibitem {tanabe}\textsc{Tanabe, H.,} \textit{Equations of Evolution,}
Monographs and Studies in Mathematics\textbf{ 6}, Pitman, London (1979).

\bibitem {tanaka}\textsc{Tanaka, H., }\textit{Stochastic Differential
Equations with Reflecting Boundary Condition in Convex Regions, }Hiroshima
Mathematical Journal \textbf{9} (1979) 163-177.

\bibitem {vancasteren}\textsc{Van Casteren, J. A.,}\textit{ Markov Processes,
Feller Semigroups and Evolution Equations,} World Scientific, Singapore (2010).

\bibitem {watson}\textsc{Watson, G. N.,} \textit{A Treatise on the Theory of
Bessel Functions, }Cambridge University Press (1996).

\bibitem {weinberger}\textsc{Weinberger, H. F.,}\textit{\ A\ First Course in
Partial Differential Equations with Complex Variables and Transform Methods},
Dover Publications Inc., New York (1995).

\bibitem {yasue}\textsc{Yasue, K., }\textit{Stochastic Calculus of
Variations,} Journal of Functional Analysis \textbf{41} (1981) 327-340.

\bibitem {zambrini}\textsc{Zambrini, J. C.,} \textit{Variational Processes and
Stochastic Versions of Mechanics,} Journal of Mathematical Physics \textbf{27
}(1986) 2307-2330.
\end{thebibliography}
\end{document}